\documentclass[a4paper,11pt]{amsart}
\usepackage{amssymb}
\numberwithin{equation}{section}
\newtheorem{theorem}{Theorem}[section]
\newtheorem{prop}[theorem]{Proposition}
\newtheorem{lem}[theorem]{Lemma}

\newtheorem{rem}[theorem]{Remark}
\theoremstyle{remark}

\newcommand{\R}{\mathbb{R}}

\newcommand{\N}{\mathbb{N}}

\author[C.~Katar]{Chaala Katar}
\address{Department of Mathematics, Faculty of Science of Gab\`es, Research Laboratory Mathematics and Applications LR17ES11; Tunisia}
\email{\sl katarchaala123@gmail.com}
\title[Asymptotic study of a global solution of super-critical ...]
{Asymptotic study of a global solution of super-critical Quasi-Geostrophic equation}

\begin{document}
\begin{abstract}
In this paper, we study the super-critical Quasi-Geostrophic equation in Gevrey-Sobolev space. We prove the local existence of $(QG)$ for any large
initial data and we give an exponential type of Blow-up to the solution. Moreover, we establish the existence global for a small initial data and we show that $\|\theta\|_{H^s_{a, \alpha^{-1}}}$ decays to zero as time goes to infinity.
Fourier analysis and standard techniques are used.
\end{abstract}
%@@@@@@@@@@@@@@@@@@@@@@@@@@@@@@@@@@@@%@@@@@@@@@@@@@@@@@@@@@@@@@@@@@@@@@@@@%@@@@@@@@@@@@@@

%\subjclass[2010]{35-XX, 35B44, 35Q30}
%\keywords{Navier-Stokes Equations; Critical spaces;  Blow up criterion}

%@@@@@@@@@@@@@@@@@@@@@@@@@@@@@@@@@@@@%@@@@@@@@@@@@@@@@@@@@@@@@@@@@@@@@@@@@%@@@@@@@@@@@@@@
\maketitle
\tableofcontents

%@@@@@@@@@@@@@@@@@@@@@@@@@@@@@@@@@@@@%@@@@@@@@@@@@@@@@@@@@@@@@@@@@@@@@@@@@%@@@@@@@@@@@@@@@

\section{Introduction}
We are concerned with the following two-dimensional quasi-geostrophic equation $(QG)$:
\begin{align*}
   (QG)\hspace {2cm}
   \begin{cases}
   \;\partial_t \theta+\kappa|D|^{2\alpha}\theta +u_{\theta}.\nabla \theta &=0\\
   \;u_{\theta} =(u^1_{\theta}, u^2_{\theta}) &=(-\partial_2|D|^{-1}\theta, \partial_1|D|^{-1}\theta)\\
   \;\theta(0, x) &=\theta^0(x).
   \end{cases}
   \end{align*}
Here  $0< \alpha < 1/2$ a real number and $\kappa > 0$ is a dissipative coefficient. The variable $\theta$ represents the potential temperature and  $u=(\partial_2|D|^{-1},\partial_1|D|^{-1} )\theta$ is the fluid velocity.\\
Our aim in this paper is to study of explosions for non-smooth solutions of $(QG)$ in finite maximal time.
 \begin{rem}
The linear system of $(QG)$ is
\begin{align*}
   (LQG)\hspace {2cm}
   \begin{cases}
   \;\partial_t \theta+|D|^{2\alpha}\theta &=0\\
  \;\theta(0, x) &=\theta^0(x) \in H^s.
   \end{cases}
   \end{align*}
Since $\theta ^0\in H^s$, the solution of the $(LQG)$ is $\theta= e^{-t|D|^{2\alpha}}\theta ^0 \in H^s_{t, (2\alpha)^{-1}}$.\\
with $H^s_{a,\sigma}$ is Gevrey Sobolev space which is defined as follow:
 for $a>0, \sigma>1 $ and $s>0$,
$$ H^s_{a,\sigma}(\mathbb{R}^d)=\{u \in L^2(\mathbb{R}^d)/ (1+|\xi|^2)^{s/2}e^{{a|\xi|}^{1/\sigma}}\in L^2(\mathbb{R}^d)\} .$$
                  equipped by the norm
                  $$\|u\|_{H^s_{a,\sigma}} =\Big(\int_{\mathbb{R}^d}(1+|\xi|^2)^s|\widehat{u}(\xi)|^2 e^{{2a|\xi|}^{1/\sigma}}d\xi\Big)^\frac{1}{2}$$
                  and the associated inner product
               $$\langle f, g\rangle_{\dot{H}^s_{a,\sigma}}= \langle e^{{a|D|}^{1/\sigma}}f,e^{{a|D|}^{1/\sigma}}g\rangle_{\dot{H}^s}$$

\end{rem}
The explosion of type exponential has been studied in the previous work of Benameur \cite{JB1}. The author use the Sobolev-Gevrey space to get better explosion result.\\
 More precisely, we use the same space in order to prove that the type of explosion is due to the chosen space not to nonlinear part of $(QG)$.\\
 %%%%%%%%%%%%%%%%%%%%%%%%%
On the other wise, we study the asymptotic behavior of the two-dimensional quasi-geostrophic equations with super critical dissipation. In literature,
the global regularity has been shown when the initial data is small in spaces $B^{2-2\alpha}_{2,1}$ \cite{ChLee}, $ H^s,s>2$ \cite{CC}, or $B^s_{2,\infty}$ with $s>2-2\alpha$ \cite{JW}..\\
 We finish by explaining why we choose $\sigma=\alpha^{-1}$ in the definition of Sobolev Gevrey  space, the reason appears in the non-linear
estimate of lemma \ref{SG}. In the proof of this lemma we have seen that in the case $\sigma=\alpha^{-1}$ there is a perfect balance between the nonlinear term and the dissipation when we have a $L^2({H^{s+\alpha}_{a,\alpha^{-1}}})$ control.\\
Let us fix $k = 1$ for the rest of the paper.\\
The following are our main theorems.
\begin{theorem}\label{TSG}
  Let $a, s, \alpha \in \mathbb{R}$ such that $a>0, s> 2$ and $ 0<\alpha<\frac{1}{2}$. Let $\theta^{0} \in H^s_{a, \alpha^{-1}}\left(\mathbb{R}^{2}\right) .$ There is a unique time $T^*\in(0,\infty]$ and a unique solution $\theta \in \mathcal{C}([0, T^*), H^s_{a, \alpha^{-1}} \left(\mathbb{R}^{2}\right))$ of $(QG)$
  Moreover, if $ T^* <\infty$, then
  \begin{align}\label{TSG.1}
  \frac{C_1}{(T^*-t)^2} \exp\big(\frac{2aC_2}{(T^*-t)^{\frac{\alpha}{s}}}\big)\leq \|\theta(t)\|^2_{\dot{H}^s_{a,\alpha^{-1}}}
  \end{align}
where $C_1$ and $C_2$ are positif constants.
\end{theorem}
For the small initial data, the global existence is given by the following theorem.
\begin{theorem}\label{TGS}
Let $\theta^0 \in  H^s_{a,\alpha^{-1}}(\mathbb{R}^2)$.\\
If $ \|\theta^0\|_{ H^s_{a,\alpha^{-1}}}<\frac{1}{2\sqrt{C}}$ then  there exists a global solution of $(QG)$ such that $$\theta \in \mathcal{C}([0,\infty), H^s_{a,\alpha^{-1}}(\mathbb{R}^2))$$
Where $C$ is positif constant.
Moreover, we have for all $t\geq0$
\begin{align}\label{RGS}
\|\theta\|_{H^s_{a,\alpha^{-1}}}^2+\int_0^t\||D|^\alpha\theta\|_{H^s_{a,\alpha^{-1}}}^2\leq \|\theta^0\|^2_{H^s_{a,\alpha^{-1}}}.
\end{align}
\end{theorem}
\begin{theorem}\label{Longtime}
Let $s>2$, $a>0$ and $ 0<\alpha<\frac{1}{2}$.\\
If $\theta \in \mathcal{C}(\mathbb{R}^+, H^s_{a,\alpha^{-1}}(\mathbb{R}^2))$ is a global solution of $(QG)$ then
$$\lim _{t\rightarrow\infty}\|\theta(t)\|^2_{H^s_{a,\alpha^{-1}}}=0.$$
\end{theorem}
\noindent The remaining part of the article is organized as follows.
\noindent In Section 2, we present some notations and we show  preliminary results which will be very useful for this paper.\\
In section 4, we prove the theorem \ref{TSG} and we have in subsection \ref{sub 6.2} a blow up result of type exponential. Then in the subsection \ref{sub 6.3}, we state that the norm of global solution in $H^s_{a, \alpha^{-1}}(\mathbb{R}^2)$ goes to zero at infinity.
\section{Notations and preliminaries results}
\subsection{Notations}
\begin{enumerate}
\item[$\bullet$]$ \mathcal{F}(f)(\xi)=\widehat{f}(\xi)=\int_{\mathbb{R}^2}\exp(-i x\xi)f(x)dx; \quad \xi=(\xi_1,\xi_2)\in \mathbb{R}^2$\\
 \item[$\bullet$] For $s\in \mathbb{R}, H^s (\mathbb{R}^2)$ denotes the usual non-homogeneous Sobolev space on $\mathbb{R}^2$
and $\langle.,.\rangle_{H^s}$ denotes the usual scalar product on $H^s(\mathbb{R}^2)$.\\
 \item[$\bullet$] For $s\in \mathbb{R}, \dot{H}^s (\mathbb{R}^2)$ denotes the usual homogeneous Sobolev space on $\mathbb{R}^2$
and $\langle.,.\rangle_{\dot{H}^s}$ denotes the usual scalar product on $\dot{H}^s(\mathbb{R}^2)$.\\
\item[$\bullet$]$ \widehat{|D|^\alpha f}(\xi)= |\xi|^\alpha\widehat{f}(\xi)$\\
\item[$\bullet$] For $\sigma\in\R$, the Fourier space is defined by $X^\sigma(\R^2)=\{f\in S'(\R^2)/\;\widehat{f}\in L^1_{loc}\;{\rm and}\;|D|^\sigma\widehat{f}\in L^1(\R^2)\}$.\\
\item[$\bullet$] $\|f\|_{X^\sigma}=\int_{\R^2}|\xi|^\sigma|\widehat{f}(\xi)|d\xi$ is a norm on $X^\sigma(\R^2)$.\\
\item[$\bullet$] $u_{\theta_1}-u_{\theta2}=u_{\theta_1-\theta_2}$.\\
\item[$\bullet$] If $(B, \|.\|)$ be a Banach space: $\mathcal{C}_b(I,B)$ = space of continuous bounded functions from an interval $I$ to $B$.
\end{enumerate}
%$\bullet$ $(X^\sigma,\|.\|_{X^\sigma})$ is a Banach space if $\sigma\leq0$. (See jamel lotfi)
\subsection{Preliminaries results}
\begin{prop}(\cite{HBAF})\label{thhb} Let $H$ be Hilbert space and $(x_n)$ be a bounded sequence of elements in $H$ such that
$$x_n\rightarrow x\;{\rm weakly\;in}\;H$$
and
$$\limsup_{n\rightarrow\infty}\|x_n\|\leq \|x\|,$$
then $$\lim_{n\rightarrow\infty}\|x_n-x\|=0.$$
\end{prop}
\begin{lem}(\cite{JB1}) Let $s\ge0$, $a>0$ and $\sigma>1$. Then, there is a constant $C=C(s)$ such that for all $f,g \in H^s_{a,\sigma}(\R^2)$, we have
$$\|fg\|_{H^s_{a,\sigma}}\leq C\Big(\|e^{\frac{a}{\sigma}|D|^{1/\sigma}}f\|_{X^0}\|g\|_{H^s_{a,\sigma}}+\|f\|_{H^s_{a,\sigma}}\|e^{\frac{a}{\sigma}|D|^{1/\sigma}}g\|_{X^0}\Big).$$
            Moreover, if $s>d/2$, we have $H^s_{a,\sigma}$ is algebra and
             $$\|fg\|_{H^s_{a,\sigma}}\leq{2c(s)}\|f\|_{H^s_{a,\sigma}}\|g\|_{H^s_{a,\sigma}}.$$
\end{lem}
\begin{lem}(\cite{JC})\label{lemjc}
          Let $0<\alpha<1/2$, there is a constant $C(\alpha)$ such that for $\theta \in H^{\sigma+\alpha}(\mathbb{R}^2)\cap H^{2-2\alpha}(\mathbb{R}^2)$ with $\sigma \geq 1$, we have
           $$|\langle u_{\theta}.\nabla\theta, \theta \rangle_{H^{\sigma}}|\leq \sigma2^\sigma C(\alpha)\|\theta\|_{\dot{H}^{2-2\alpha}} \|\theta\|^2_{\dot{H}^{\sigma+\alpha}}.$$
 \end{lem}
\begin{lem}\label{SG}
   Let $s>2$, $a>0$, $\alpha\in(0,1)$ and for every $\theta, \omega \in H^{s+\alpha}_{a,\alpha^{-1}}(\mathbb{R}^2)$, there is a constant $C=  C_{s,a, \alpha}$ such that
\begin{equation}\label{eql11}
|\langle u_{\theta}.\nabla\omega, \omega \rangle_{H^s_{a,\alpha^{-1}}}|\leq C\Big(\|e^{a\alpha|D|^\alpha}\theta\|_{X^1}\||D|^\alpha\omega\|_{H^s_{a,\alpha^{-1}}} +\||D|^\alpha\theta\|_{H^s_{a,\alpha^{-1}}}\|e^{a\alpha|D|^\alpha}\omega\|_{X^1}\Big)\|\omega\|_{H^s_{a,\alpha^{-1}}},
\end{equation}
\begin{equation}\label{eql13}
|\langle u_{\theta}.\nabla\theta, \theta \rangle_{H^s_{a,\alpha^{-1}}}|\leq C\|e^{a\alpha|D|^\alpha}\theta)\|_{X^1}\||D|^\alpha\theta\|_{H^s_{a,\alpha^{-1}}}\|\theta\|_{H^s_{a,\alpha^{-1}}},
\end{equation}
\begin{equation}\label{eql14}
|\langle u_{\theta}.\nabla\theta, \theta \rangle_{H^s_{a,\alpha^{-1}}}|\leq C\||D|^\alpha\theta\|_{H^s_{a,\alpha^{-1}}}\|\theta\|_{H^s_{a,\alpha^{-1}}}^2,
\end{equation}
\begin{equation}\label{eql12}
|\langle u_{\theta}.\nabla\omega, \omega \rangle_{H^s_{a,\alpha^{-1}}}|\leq C\|\theta\|_{H^s_{a,\alpha^{-1}}} \||D|^\alpha\omega\|_{H^s_{a,\alpha^{-1}}}\|\omega\|_{H^s_{a,\alpha^{-1}}}.
\end{equation}
\end{lem}
{\bf Proof.} Clearly equation (\ref{eql13}) is a particular case of (\ref{eql11}). We start by proving the first equation.\\
$\bullet$ {Proof of (\ref{eql11}):} Using the fact $\langle u_\theta.\nabla\omega, \omega\rangle_{L^2}=0$ and $\langle f, g\rangle_{H^s_{a,\sigma}}= \langle f, g\rangle_{L^2}+\langle f, g\rangle_{\dot{H}^s_{a,\sigma}},$ we get$$< u_{\theta}.\nabla\omega, \omega>_{H^s_{a,\sigma}}=<u_{\theta}.\nabla\omega, \omega>_{\dot H^s_{a,\sigma}}.$$
Using again $\langle u_\theta.\nabla |D|^se^{a|D|^\alpha}\omega, |D|^se^{a|D|^\alpha}\omega\rangle_{L^2}=0$ and Cauchy-Schwartz inequality, we obtain
$$\begin{array}{lcl}
|\langle u_{\theta}.\nabla\omega, \omega\rangle_{\dot{H}^s_{a,\alpha^{-1}}}|&\leq&
\displaystyle\int_\xi\int||\xi|^s e^{a|\xi|^{\alpha}}-|\eta|^s e^{a|\eta|^{\alpha}}|.|\widehat{\theta}(\xi-\eta)||\widehat{\nabla\omega}(\eta)||\xi|^s e^{a|\xi|^\alpha}|\widehat{\omega}(-\xi)|d\xi\\
&\leq&\displaystyle\Big(\int_\xi(\int||\xi|^s e^{a|\xi|^{\alpha}}-|\eta|^s e^{a|\eta|^{\alpha}}|.|\widehat{\theta}(\xi-\eta)||\widehat{\nabla\omega}(\eta)|d\eta)^2d\xi\Big)^{1/2}\|\omega\|_{\dot H^s_{a,\alpha^{-1}}}.
\end{array}$$
For $\xi,\;\eta\in\R^2$, there is $z\in[\min(|\xi|,|\eta|),\max(|\xi|,|\eta|)]$ such that
$$|\xi|^s e^{a|\xi|^\alpha}-|\eta|^s e^{a|\eta|^\alpha}=(|\xi|-|\eta|)\Big(sz^{s-1}+a\alpha z^{s+\alpha-1}\Big)e^{a z^\alpha}.$$
Using the fact $z\leq \max(|\xi|,|\eta|)\leq |\xi-\eta|+|\eta|\leq 2\max(|\xi-\eta|,|\eta|)$ and the fact (see \cite{JB1}) $$|\xi|^\alpha\leq \max(|\xi-\eta|,|\eta|)^\alpha+\alpha\min(|\xi-\eta|,|\eta|)^\alpha,$$ we get
$$sz^{s-1}+a\alpha z^{s+\alpha-1}\leq 2^{s+\alpha-1}(s+a\alpha)\Big(\max(|\xi-\eta|,|\eta|)^{s-1}+\max(|\xi-\eta|,|\eta|)^{s+\alpha-1}\Big)$$
$$\begin{array}{lcl}e^{az^\alpha}&\leq&\left\{\begin{array}{l}
e^{a|\eta|^\alpha},\;{\rm if}\;|\xi|<|\eta|\\
e^{a|\xi|^\alpha},\;{\rm if}\;|\xi|>|\eta|\\
\end{array}\right.\\\\
&\leq&\left\{\begin{array}{l}
e^{a|\eta|^\alpha},\;{\rm if}\;|\xi|<|\eta|\\
e^{a\max(|\xi-\eta|,|\eta|)^\alpha+a\alpha\min(|\xi-\eta|,|\eta|)^\alpha},\;{\rm if}\;|\xi|>|\eta|\\
\end{array}\right.\\\\
&\leq&e^{a\max(|\xi-\eta|,|\eta|)^\alpha+a\alpha\min(|\xi-\eta|,|\eta|)^\alpha}.\end{array}
$$
Then
$$\Big||\xi|^s e^{a|\xi|^\alpha}-|\eta|^s e^{a|\eta|^\alpha}\Big|\leq
\left\{\begin{array}{l}
C|\xi-\eta|\Big(|\xi-\eta|^{s-1}+|\xi-\eta|^{s+\alpha-1}\Big)e^{a|\xi-\eta|^\alpha}e^{a\alpha|\eta|^\alpha},\;{\rm if}\;|\xi-\eta|>|\eta|\\
C|\xi-\eta|\Big(|\eta|^{s-1}+|\eta|^{s+\alpha-1}\Big)e^{a|\eta|^\alpha}e^{a\alpha|\xi-\eta|^\alpha},\;{\rm if}\;|\xi-\eta|<|\eta|\\
\end{array}\right.$$
with $C=2^{s+\alpha}(s+\alpha)$. Therefore
$$|\langle u_{\theta}.\nabla\omega, \omega\rangle_{\dot{H}^s_{a,\alpha^{-1}}}|\leq C\Big(\sum_{k=1}^4I_k\Big)\|\omega\|_{\dot H^s_{a,\alpha^{-1}}},$$
with
$$\begin{array}{lcl}
I_1&=&\displaystyle\Big(\int_\xi(\int_\eta|\xi-\eta|^s e^{a|\xi-\eta|^{\alpha}} |\widehat{\theta}(\xi-\eta)|e^{a\alpha|\eta|^{\alpha}}|\eta|.|\widehat{\omega}(\eta)|d\eta)^2d\xi\Big)^{1/2}:=\|F_1*G_1\|_{L^2}\\
&&F_1=|\xi|^s e^{a|\xi|^{\alpha}} |\widehat{\theta}(\xi)|,\;\;G_1=|\xi| e^{a\alpha|\xi|^{\alpha}} |\widehat{\omega}(\xi)|\\
I_2&=&\displaystyle\Big(\int_\xi(\int_\eta|\xi-\eta|^{s+\alpha} e^{a|\xi-\eta|^{\alpha}} |\widehat{\theta}(\xi-\eta)|e^{a\alpha|\eta|^{\alpha}}|\eta|.|\widehat{\omega}(\eta)|d\eta)^2d\xi\Big)^{1/2}:=\|F_2*G_2\|_{L^2}\\
&&F_2=|\xi|^{s+\alpha} e^{a|\xi|^{\alpha}} |\widehat{\theta}(\xi)|,\;\;G_2=|\xi| e^{a\alpha|\xi|^{\alpha}} |\widehat{\omega}(\xi)|\\
I_3&=&\displaystyle\Big(\int_\xi(\int_\eta|\xi-\eta| e^{a\alpha|\xi-\eta|^{\alpha}} |\widehat{\theta}(\xi-\eta)|e^{a|\eta|^{\alpha}}|\eta|^s|\widehat{\omega}(\eta)|d\eta)^2d\xi\Big)^{1/2}:=\|F_3*G_4\|_{L^2}\\
&&F_3=|\xi|e^{a\alpha|\xi|^{\alpha}} |\widehat{\theta}(\xi)|,\;\;G_3=|\xi|^s e^{a|\xi|^{\alpha}} |\widehat{\omega}(\xi)|\\
I_4&=&\displaystyle\Big(\int_\xi(\int_\eta|\xi-\eta| e^{a\alpha|\xi-\eta|^{\alpha}} |\widehat{\theta}(\xi-\eta)|.|\eta|^{s+\alpha}e^{a|\eta|^{\alpha}}|\widehat{\omega}(\eta)|d\eta)^2d\xi\Big)^{1/2}:=\|F_4*G_4\|_{L^2}\\
&&F_4=|\xi| e^{a\alpha|\xi|^{\alpha}} |\widehat{\theta}(\xi)|,\;\;G_4=|\xi|^{s+\alpha} e^{a|\xi|^{\alpha}} |\widehat{\omega}(\xi)|.
\end{array}$$
Young inequality implies
$$|\langle u_{\theta}.\nabla\omega, \omega\rangle_{\dot{H}^s_{a,\sigma}}|\leq C\Big(\|F_1\|_{L^2}\|G_1\|_{L^1}+\|F_2\|_{L^2}\|G_2\|_{L^1}+\|F_3\|_{L^1}\|G_1\|_{L^2}+\|F_4\|_{L^1}\|G_4\|_{L^2}\Big)\|\omega\|_{\dot H^s_{a,\alpha^{-1}}}.$$
As $0\leq \alpha\leq s$, then
$$|\xi|^s\leq \left\{\begin{array}{l}
|\xi|^\alpha\;{\rm if}\;|\xi|<1\\
|\xi|^{s+\alpha}\;{\rm if}\;|\xi|>1
\end{array}\right.\leq (1+|\xi|^2)^{\frac{s}{2}}|\xi|^\alpha.$$
By the above inequality we obtain
$$\begin{array}{lcl}
\|F_1\|_{L^2}&=&\|\theta\|_{\dot H^s_{a,\alpha^{-1}}}\leq \||D|^\alpha\theta\|_{H^s_{a,\alpha^{-1}}}\\
\|G_3\|_{L^2}&=&\|\omega\|_{\dot H^s_{a,\alpha^{-1}}}\leq \||D|^\alpha\omega\|_{H^s_{a,\alpha^{-1}}}.
\end{array}$$
which imply the desired result.\\
\noindent$\bullet$ {Proof of (\ref{eql14}):} By using equation (\ref{eql13}), it suffices  to prove $\|e^{a\alpha|D|^\alpha}\theta)\|_{X^1}\leq C\|\theta)\|_{H^s_{a,\alpha^{-1}}}$. For this, we write
$$\begin{array}{lcl}
\|e^{a\alpha|D|^\alpha}\theta\|_{X^1}&=&\displaystyle\int_{\R^2}|\xi|e^{a\alpha|\xi|^\alpha}|\widehat{\theta}(\xi)|d\xi\\
&=&\displaystyle\int_{\R^2}\frac{|\xi|e^{a\alpha|\xi|^\alpha}}{(1+|\xi|^2)^{s/2}e^{a|\xi|^\alpha}}(1+|\xi|^2)^{s/2}e^{a|\xi|^\alpha}|\widehat{\theta}(\xi)|d\xi\\
&\leq&\displaystyle\Big(\int_{\R^2}\frac{|\xi|e^{a\alpha|\xi|^\alpha}}{(1+|\xi|^2)^{s/2}e^{a|\xi|^\alpha}}d\xi\Big)^{1/2}\|\theta\|_{\dot H^s_{a,\alpha^{-1}}},
\end{array}$$
which gives (\ref{eql14}).\\
\noindent$\bullet$ {Proof of (\ref{eql12}):} We have
$$|\langle u_{\theta}.\nabla\omega, \omega\rangle_{\dot{H}^s_{a,\sigma}}|\leq
\displaystyle\int_\xi\int||\xi|^s e^{a|\xi|^{\alpha}}-|\eta|^s e^{a|\eta|^{\alpha}}|.|\widehat{\theta}(\xi-\eta)||\widehat{\nabla\omega}(\eta)||\xi|^s e^{a|\xi|^\alpha}|\widehat{\omega}(-\xi)|d\xi.$$
By using the same steps of the above proof, we get: For $\xi,\;\eta\in\R^2$,
$$\Big||\xi|^s e^{a|\xi|^\alpha}-|\eta|^s e^{a|\eta|^\alpha}\Big|\leq
\left\{\begin{array}{l}
C|\xi-\eta|\Big(|\xi-\eta|^{s-1}+|\xi-\eta|^{s+\alpha-1}\Big)e^{a|\xi-\eta|^\alpha}e^{a\alpha|\eta|^\alpha},\;{\rm if}\;|\xi-\eta|>|\eta|\\
C|\xi-\eta|\Big(|\eta|^{s-1}+|\eta|^{s+\alpha-1}\Big)e^{a|\eta|^\alpha}e^{a\alpha|\xi-\eta|^\alpha},\;{\rm if}\;|\xi-\eta|<|\eta|.
\end{array}\right.$$
Using the fact $|\xi-\eta|^{s+\alpha-1}\leq c|\xi-\eta|^{s-1}(|\xi|^\alpha+|\eta|^\alpha)$ we get
$$\Big||\xi|^s e^{a|\xi|^\alpha}-|\eta|^s e^{a|\eta|^\alpha}\Big|\leq
\left\{\begin{array}{l}
C|\xi-\eta|\Big(|\xi-\eta|^{s-1}+|\xi-\eta|^{s-1}(|\xi|^\alpha+|\eta|^\alpha)\Big)e^{a|\xi-\eta|^\alpha}e^{a\alpha|\eta|^\alpha},\;{\rm if}\;|\xi-\eta|>|\eta|\\\\
C|\xi-\eta|\Big(|\eta|^{s-1}+|\eta|^{s+\alpha-1}\Big)e^{a|\eta|^\alpha}e^{a\alpha|\xi-\eta|^\alpha},\;{\rm if}\;|\xi-\eta|<|\eta|.
\end{array}\right.$$
Therefore
$$|\langle u_{\theta}.\nabla\omega, \omega\rangle_{\dot{H}^s_{a,\sigma}}|\leq C\sum_{k=1}^5J_k,$$
with
$$\begin{array}{lcl}
J_1&=&\displaystyle\int_\xi(\int_\eta|\xi-\eta|^s e^{a|\xi-\eta|^{\alpha}}|\widehat{\theta}(\xi-\eta)|e^{a\alpha|\eta|^{\alpha}}|\eta|.|\widehat{\omega}(\eta)|d\eta)|\xi|^s e^{a|\xi|^\alpha}|\widehat{\omega}(-\xi)|d\xi\\
J_2&=&\displaystyle\int_\xi(\int_\eta|\xi-\eta|^s e^{a|\xi-\eta|^{\alpha}}|\widehat{\theta}(\xi-\eta)|e^{a\alpha|\eta|^{\alpha}}|\eta|.|\widehat{\omega}(\eta)|d\eta)|\xi|^{s+\alpha} e^{a|\xi|^\alpha}|\widehat{\omega}(-\xi)|d\xi\\
J_3&=&\displaystyle\int_\xi(\int_\eta|\xi-\eta|^s e^{a|\xi-\eta|^{\alpha}}|\widehat{\theta}(\xi-\eta)|e^{a\alpha|\eta|^{\alpha}}|\eta|^{1+\alpha}.|\widehat{\omega}(\eta)|d\eta)|\xi|^s e^{a|\xi|^\alpha}|\widehat{\omega}(-\xi)|d\xi\\
J_4&=&\displaystyle\int_\xi(\int_\eta|\xi-\eta| e^{a\alpha|\xi-\eta|^{\alpha}}|\widehat{\theta}(\xi-\eta)|e^{a|\eta|^{\alpha}}|\eta|^s|\widehat{\omega}(\eta)|d\eta)|\xi|^s e^{a|\xi|^\alpha}|\widehat{\omega}(-\xi)|d\xi\\
J_5&=&\displaystyle\int_\xi(\int_\eta|\xi-\eta| e^{a\alpha|\xi-\eta|^{\alpha}}|\widehat{\theta}(\xi-\eta)|e^{a|\eta|^{\alpha}}|\eta|^{s+\alpha}|\widehat{\omega}(\eta)|d\eta)|\xi|^s e^{a|\xi|^\alpha}|\widehat{\omega}(-\xi)|d\xi.
\end{array}$$
Combining Cauchy-Schwartz and Young inequalities, we get
$$\begin{array}{lcl}
J_1&\leq&\displaystyle\|\theta\|_{\dot H^s_{a,\alpha^{-1}}}\|e^{a\alpha|D|^{\alpha}}\omega\|_{X^1}\|\omega\|_{\dot H^s_{a,\alpha^{-1}}}\\
J_2&\leq&\displaystyle\|\theta\|_{\dot H^s_{a,\alpha^{-1}}}\|e^{a\alpha|D|^{\alpha}}\omega\|_{X^1}\|\omega\|_{\dot H^{s+\alpha}_{a,\alpha^{-1}}}\\
J_3&\leq&\displaystyle\|\theta\|_{\dot H^s_{a,\alpha^{-1}}}\|e^{a\alpha|D|^{\alpha}}\omega\|_{X^{1+\alpha}}\|\omega\|_{\dot H^s_{a,\alpha^{-1}}}\\
J_4&\leq&\displaystyle\|e^{a\alpha|D|^{\alpha}}\theta\|_{X^1}\|\omega\|_{\dot H^s_{a,\alpha^{-1}}}\|\omega\|_{\dot H^s_{a,\alpha^{-1}}}\\
J_5&\leq&\displaystyle\|e^{a\alpha|D|^{\alpha}}\theta\|_{X^1}\|\omega\|_{\dot H^{s+\alpha}_{a,\alpha^{-1}}}\|\omega\|_{\dot H^s_{a,\alpha^{-1}}}.
\end{array}$$
Using the fact $|\xi|^s\leq (1+|\xi|^s)|\xi|^\alpha, \quad (0<\alpha<s)$ and
$$\int_{\R^2}\frac{|\xi|^2e^{2a\alpha|\xi|^{\alpha}}}{(1+|\xi|^2)^se^{2a|\xi|^{\alpha}}}d\xi,\;
\int_{\R^2}\frac{|\xi|^{2+2\alpha}e^{2a\alpha|\xi|^{\alpha}}}{(1+|\xi|^2)^se^{2a|\xi|^{\alpha}}}d\xi<\infty$$
we obtain the result.
%%%%%%%%%%%%%%%%%%%%%%%%%%%%%%%
\begin{lem}
Let $0<\alpha<\frac{1}{2}$, there is a constant $C(s)$ such that for $\theta \in H^{s+\alpha}_{a,\alpha^{-1}}(\mathbb{R}^2)$ with $s \geq 2$, we have
We have
\begin{equation}\label{eql15}
|\langle u_{\theta}.\nabla\theta, \theta \rangle_{H^s_{a,\alpha^{-1}}}|\leq C\|\theta\|_{H^s_{a,\alpha^{-1}}}\||D|^\alpha\theta\|^2_{H^s_{a,\alpha^{-1}}},
\end{equation}
\end{lem}
{\bf proof.}
Using $\langle u_\theta.\nabla |D|^se^{a|D|^\alpha}\theta, |D|^se^{a|D|^\alpha}\theta\rangle_{L^2}=0$, we obtain
$$\begin{array}{lcl}
|\langle u_{\theta}.\nabla\theta, \theta\rangle_{\dot{H}^s_{a,\alpha^{-1}}}|&\leq&
\displaystyle\int_\xi\int_\eta||\xi|^s e^{a|\xi|^{\alpha}}-|\eta|^s e^{a|\eta|^{\alpha}}|.|\widehat{\theta}(\xi-\eta)||\widehat{\nabla\theta}(\eta)||\xi|^s e^{a|\xi|^\alpha}|\widehat{\alpha}(-\xi)|d\eta d\xi.
\end{array}$$
The inequality $e^{a|\xi|^{\alpha}}\leq e^{a|\xi-\eta|^{\alpha}}e^{a|\eta|^{\alpha}}$ implies
$$\Big||\xi|^s e^{a|\xi|^\alpha}-|\eta|^s e^{a|\eta|^\alpha}\Big|\leq
\left\{\begin{array}{l}
C|\xi-\eta|\Big(|\xi-\eta|^{s-1}+|\xi-\eta|^{s+\alpha-1}\Big)e^{a|\xi-\eta|^\alpha}e^{a|\eta|^\alpha},\;{\rm if}\;|\xi-\eta|>|\eta|\\
C|\xi-\eta|\Big(|\eta|^{s-1}+|\eta|^{s+\alpha-1}\Big)e^{a|\eta|^\alpha}e^{a|\xi-\eta|^\alpha},\;{\rm if}\;|\xi-\eta|<|\eta|\\
\end{array}\right.$$
Therefore
          $$|\langle u_{\theta}.\nabla\theta, \theta \rangle_{H^s_{a,\alpha^{-1}}}|
         \leq C \sum_{k=1}^4I_k,$$
with
         $$\begin{array}{lcl}
I_1&\leq&\displaystyle\int_{\xi}(\int_{\eta}|\xi-\eta|^se^{a|\xi-\eta|^{\alpha}}|\widehat{\theta}(\xi-\eta)|.|\eta|e^{a|\eta|^{\alpha}}|\widehat{\theta}(\eta)|d\eta)|\xi|^se^{a|\xi|^{\alpha}}|\widehat{\theta}(\xi)|d\xi\\
I_2&\leq&\displaystyle\int_{\xi}(\int_{\eta}|\xi-\eta|^{s+\alpha}e^{a|\xi-\eta|^{\alpha}}|\widehat{\theta}(\xi-\eta)|.|\eta|e^{a|\eta|^{\alpha}}|\widehat{\theta}(\eta)|d\eta)|\xi|^s|\widehat{\theta}(\xi)|d\xi\\
I_3&\leq&\displaystyle\int_{\xi}(\int_{\eta}|\xi-\eta|e^{a|\xi-\eta|^{\alpha}}|\widehat{\theta}(\xi-\eta)|.|\eta|^se^{a|\eta|^{\alpha}}|\widehat{\theta}(\eta)|d\eta)|\xi|^se^{a|\xi|^{\alpha}}|\widehat{\theta}(\xi)|d\xi\\
I_4&\leq&\displaystyle\int_{\xi}(\int_{\eta}|\xi-\eta|e^{a|\xi-\eta|^{\alpha}}|\widehat{\theta}(\xi-\eta)|.|\eta|^{s+\alpha}e^{a|\eta|^{\alpha}}|\widehat{\theta}(\eta)|d\eta)|\xi|^s|\widehat{\theta}(\xi)|d\xi.\\
\end{array}$$
Then \begin{align*}
 |\langle u_{\theta}.\nabla\theta, \theta \rangle_{H^s_{a,\alpha^{-1}}}|\leq C\big(\|f_1g_1\|_{\dot H^{-\alpha}}+\|f_2g_2\|_{\dot H^{-\alpha}}\big) \|\theta\|_{\dot{H}^{s+\alpha}}
 \end{align*}
 with
 $$\left\{\begin{array}{lcl}
\widehat{f_1}&=&|\xi|e^{a|\xi|^{\alpha}}|\widehat{\theta}(\xi)|,\\
\widehat{g_1}&=&|\xi|^se^{a|\xi|^{\alpha}}|\widehat{\theta}(\xi)|,\\
\widehat{f_2}&=&|\xi|e^{a|\xi|^{\alpha}}|\widehat{\theta}(\xi)|,\\
\widehat{g_2}&=&|\xi|^{s+\alpha}e^{a|\xi|^{\alpha}}|\widehat{\theta}(\xi)|.
\end{array}\right.$$
Using the product laws in homogenous Sobolev space for $f_1, g_1$, with $ s_1+s_2=1-\alpha>0$ and
\begin{align*}
\begin{cases}
&s_1=1-2\alpha<1\\
&s_2=\alpha<1
\end{cases}
\end{align*}
and for $f_2, g_2$, with $ s_1+s_2=1-\alpha>0$ and
\begin{align*}
\begin{cases}
&s_1=1-\alpha<1\\
&s_2=0<1
\end{cases}
\end{align*}
Combining the above results, we get
\begin{align*}
|\langle u_{\theta}.\nabla\theta, \theta \rangle_{H^s_{a,\alpha^{-1}}}|&\leq 2C\big(\|f_1\|_{\dot{H}^{1-2\alpha}} \|g_1\|_{\dot{H}^{\alpha}}+\|f_2\|_{\dot{H}^{1-\alpha}} \|g_2\|_{\dot{H}^{0}}\big)\|\theta\|_{\dot{H}^{s+\alpha}}\\
&\leq\big(\|\theta\|_{\dot H^{2-2\alpha}_{a,\alpha^{-1}}}\|\theta\|_{\dot H^{s+\alpha}_{a,\alpha^{-1}}}
+\|\theta\|_{\dot H^{2-\alpha}_{a,\alpha^{-1}}} \|\theta\|_{\dot H^{s+\alpha}_{a,\alpha^{-1}}}
\big)\|\theta\|_{\dot H^{s+\alpha}_{a,\alpha^{-1}}}.\\
&\leq\big(\|\theta\|_{ H^{2-2\alpha}_{a,\alpha^{-1}}}\|\theta\|_{H^{s+\alpha}_{a,\alpha^{-1}}}
+\|\theta\|_{H^{2-\alpha}_{a,\alpha^{-1}}} \|\theta\|_{ H^{s+\alpha}_{a,\alpha^{-1}}}
\big)\|\theta\|_{ H^{s+\alpha}_{a,\alpha^{-1}}}.\\
\end{align*}
Now, using the fact
$2-2\alpha\leq 2<s$ and $2-\alpha\leq 2<s$ to obtain the desired result.

  %%%%%%%%%%%%%%%%%%%%%%%%%%%%%%%%%%%%%%%%
\begin{rem}
For $a>0, s>0$ and $\sigma>1$, we have
               $$\overline{e^{-a|D|^{1/\sigma}} \mathcal{S}(\mathbb{R}^2)}=H^s_{a,\sigma}(\mathbb{R}^2).$$
               \end{rem}
Indeed, it suffices to write $$\overline{\mathcal{S}(\mathbb{R}^2)}=H^s(\mathbb{R}^2) \quad \forall s>0.$$
\subsection{Local well-posedness}\label{sub 6.1}
This subsection is devoted to the proof of the local well-posedness of the system $(QG)$, which we give the proof the following theorem.
\begin{theorem}
  Let $a, s, \alpha \in \mathbb{R}$ such that $a>0, s> 2$ and $0<\alpha<\frac{1}{2}$. Let $\theta^{0} \in H^s_{a, \alpha^{-1}}\left(\mathbb{R}^{2}\right) .$ There is a time $T>0$ and a unique solution $\theta \in \mathcal{C}\left([0, T], H^s_{a, \alpha^{-1}} \left(\mathbb{R}^{2}\right)\right)$ of $(QG)$.
\end{theorem}
{\bf Proof.}
Proving this theorem requires four steps:\\
- First, we apply the method fixed point theorem so as to solve a approximating system $(QG)_k$ of $(QG)$.\\
-  Second, we prove that for uniform T sufficiently small, the sequence $(\theta_k)$ is bounded in $\mathcal{C}\left([0, T] ; H^{s}_{a, \alpha^{-1}} \right)$.\\
- Third, we establish that for T sufficiently small, $(\theta_k)$ is a Cauchy sequence in $L^{\infty}([0, T] ; H^{s-1}_{a, \alpha^{-1}} )$.\\
- Finally, we check that the limit of this sequence is a solution of $(QG)$ and that it belongs to $\mathcal{C}([0,T];H^{s}_{a, \alpha^{-1}})$.\\
\textbf{Step1: Construction of an approximate solution sequence }\\
We will use the Kato method \cite{K1} to construct the approximate solutions.\\
We introduce the following approximate system of $(QG)$, for $ k\in \N$,
\begin{align*}
   (QG)_k \hspace {2cm}
   \begin{cases}
   \;\partial_t \theta-\frac{1}{k}\Delta \theta+(-\Delta)^{\alpha}\theta +u_{\theta}.\nabla \theta &=0\\
   \;\theta_k(0, x) &=\theta^0(x).
   \end{cases}
   \end{align*}
The existence of solution of $(QG)_k$ based on fixed point theorem.\\
$(QG)_k$ has following integral form
$$\theta _k = e^{t(\frac{1}{k} \Delta- (-\Delta)^\alpha)}\theta ^0-\int_0^te^{(t-\tau)(\frac{1}{k} \Delta- (-\Delta)^\alpha)}(u_{\theta}.\nabla \theta)d\tau $$
Let $$X_{T}=
\left\{u \in C\left([0, T] ;H_{a, \alpha^{-1}}^{s}\left(\mathbb{R}^{2}\right)\right) ;\|u\|_{L^{\infty}_T(H_{a, \alpha^{-1}}^{s}).} \leq 2\|\theta^0\|_{H^s_{a, \alpha^{-1}}}\right\}
$$
 and we consider $$
\begin{aligned} F &: X_{ T} \rightarrow C([0, T] ;H_{a, \alpha^{-1}}^{s}(\mathbb{R}^{2})) \\ & \theta \mapsto e^{t(\frac{1}{k} \Delta- (-\Delta)^\alpha)}\theta ^0-\int_0^te^{(t-\tau)(\frac{1}{k} \Delta- (-\Delta)^\alpha)}(u_{\theta}.\nabla \theta)d\tau  \end{aligned}
$$
$\bullet $ Firstly, we prove $F (X_{T})\subset X_{T}$, we have
$$\|e^{t(\frac{1}{k} \Delta- (-\Delta)^\alpha)}\theta ^0\|_{H^{s}_{a, \alpha^{-1}}}\leq \|\theta ^0\|_{H^{s}_{a, \alpha^{-1}}}.$$
To estimate the non linear part, we can write:
$$\|\int_0^te^{(t-\tau)(\frac{1}{k} \Delta- (-\Delta)^\alpha)}(u_{\theta}.\nabla \theta)d\tau\|_{H^{s}_{a, \alpha^{-1}}}\leq\int_0^t\|e^{(t-\tau)(\frac{1}{k} \Delta- (-\Delta)^\alpha)}(u_{\theta}.\nabla \theta)\|_{H^{s}_{a, \alpha^{-1}}}d\tau.$$
We have
\begin{align*}
\|e^{(t-\tau)(\frac{1}{k} \Delta- (-\Delta)^\alpha)}(u_{\theta}.\nabla \theta)\|^2_{H^{s}_{a, \alpha^{-1}}}&\leq\int\left(1+|\xi|^{2}\right)^{s} e^{2 a|\xi|^{\alpha}}\left|\mathcal{F}\left(e^{(t-\tau)(\frac{1}{k} \Delta- (-\Delta)^\alpha)} \operatorname{div}(\theta u_\theta)\right)(\tau, \xi)\right|^{2} d \xi\\
&\leq\int\left(1+|\xi|^{2}\right)^{s} e^{2 a|\xi|^{\alpha}}e^{-2(t-\tau)(\frac{1}{k} |\xi|^2+|\xi|^{2\alpha})} |\xi|^2|\widehat{\theta u_\theta}|^2d\xi\\
&\leq\int\left(1+|\xi|^{2}\right)^{s} e^{2 a|\xi|^{\alpha}}e^{-2(t-\tau)\frac{1}{k} |\xi|^2} 2(t-\tau)\frac{1}{k} |\xi|^2\frac{k}{2(t-\tau)}|\widehat{\theta u_\theta}|^2d\xi\\
&\leq \frac{k}{2(t-\tau)}\|\theta u_\theta\|^2_{H^{s}_{a, \alpha^{-1}}}\\
&\leq C(a,s,\alpha) \frac{k}{2(t-\tau)}\|u_\theta\|^2_{H^{s}_{a, \alpha^{-1}}}\|\theta\|^2_{H^{s}_{a, \alpha^{-1}}},\quad \ (s>1),
\end{align*}
then
\begin{align*}
\|\int_0^te^{(t-\tau)(\frac{1}{k} \Delta- (-\Delta)^\alpha)}(u_{\theta}.\nabla \theta)d\tau\|_{H^{s}_{a, \alpha^{-1}}}&\leq\int_0^t C_s\left(\frac{k}{2(t-\tau)}\right)^{\frac{1}{2}}\|\theta\|^2_{H^{s}_{a, \alpha^{-1}}}d\tau\\
&\leq C(a,s,\alpha) \sqrt{k}\sqrt{T}\|\theta\|^2_{L^{\infty}_T( H^{s}_{a, \alpha^{-1}})}.
\end{align*}
Then, for $T<0$ such that
$$ \|\theta^0\|^2_{H^{s}_{a, \alpha^{-1}}}+  C_s\sqrt{k}\sqrt{T}\|\theta\|^2_{L^{\infty}_T( H^{s}_{a, \alpha^{-1}})}\leq 2\|\theta^0\|_{H^s_{a, \alpha^{-1}}}.$$
 We get $ F(X_{T})\subset X_{T}.$\\
$\bullet$ We show that $F$ is contracting with an additional condition over $T$.\\
We have
\begin{align*}
\|F(\theta_1)(t)-F(\theta_2)(t)\|_{H^{s}_{a, \alpha^{-1}}}&\leq \|\int_0^T e^{(t-\tau)(\frac{1}{k} \Delta- (-\Delta)^\alpha)}(u_{\theta_1}.\nabla \theta_1-u_{\theta_2}.\nabla \theta_2)d\tau\|_{H^{s}_{a, \alpha^{-1}}}\\
&\leq \|\int_0^T e^{(t-\tau)(\frac{1}{k} \Delta- (-\Delta)^\alpha)}\left(u_{\theta_1}.\nabla (\theta_1-\theta_2)+(u_{\theta_2}-u_{\theta_1}).\nabla \theta_2\right)d\tau\|_{H^{s}_{a, \alpha^{-1}}}\\
&\leq K_1+K_2,
\end{align*}
with \begin{align*}
K_1=\|\int_0^T e^{(t-\tau)(\frac{1}{k} \Delta- (-\Delta)^\alpha)}u_{\theta_1}.\nabla (\theta_1-\theta_2)d\tau\|_{H^{s}_{a, \alpha^{-1}}},\\
K_2=\|\int_0^T e^{(t-\tau)(\frac{1}{k} \Delta- (-\Delta)^\alpha)}(u_{\theta_2}-u_{\theta_1}).\nabla \theta_2d\tau\|_{H^{s}_{a, \alpha^{-1}}}.
\end{align*}
Since $\theta_1, \theta_2 \in X_{R,T}$ and $\|u_\theta\|_{H^{s}_{a, \alpha^{-1}}}=\|\theta\|_{H^{s}_{a, \alpha^{-1}}}$, we obtain \begin{align*}
\begin{cases}
K_1\leq C_s\sqrt{k}\sqrt{T}\|\theta_1-\theta_2\|_{H^{s}_{a, \alpha^{-1}}}\\
K_2\leq C_s\sqrt{k}\sqrt{T}\|\theta_1-\theta_2\|_{H^{s}_{a, \alpha^{-1}}}
\end{cases}
\end{align*}
Then $$\|F(\theta_1)(t)-F(\theta_2)(t)\|_{H^{s}_{a, \alpha^{-1}}}\leq 2C_s\sqrt{k}\sqrt{T}\|\theta_1-\theta_2\|_{L^{\infty}_T( H^{s}_{a, \alpha^{-1}})}$$
We choose $T>0$ such that $$2C_s\sqrt{k}\sqrt{T}<1$$ in order to obtain the contracting case of $F$.\\
For these choices,  fixed point theorem guarantee the existence unique solution $\theta_k$ of $(QG)_k$ in $\mathcal{C}([0,T], H^{s}_{a, \alpha^{-1}}).$\\
\noindent \textbf{Step2: Energy estimates and uniform time}\\
Let $\theta _k \in \mathcal{C}([0, T^*_k),H^{s}_{a, \alpha^{-1}})$ the maximal solution of $(QG)_k$, given by the first step with $T^*_k \in (0,\infty]$.
Taking the scalar product in $H^{s}_{a, \alpha^{-1}}$ and taking into account of (2.3), we obtain
\begin{align*}
\frac{1}{2}\frac{d}{dt}\|\theta_k\|^2_{H^s_{a,\alpha^{-1}}}+\frac{1}{k}\||D|\theta_k\|^2_{H^s_{a,\alpha^{-1}}}+\||D|^\alpha \theta_k\|^2_{H^s_{a,\alpha^{-1}}}&\leq |\langle u_{\theta_k}.\nabla\theta_k, \theta_k\rangle|_{H^s_{a,\alpha^{-1}}}\\
&\leq C\|\theta_k\|^2_{H^s_{a,\alpha^{-1}}} \ \ \||D|^\alpha\theta_k\|_{H^s_{a,\alpha^{-1}}}.
\end{align*}
 By the convex inequality $ab\leq \frac{a^2}{2}+\frac{b^2}{2}$, we have
$$\frac{1}{2}\frac{d}{dt}\|\theta_k\|^2_{H^s_{a,\alpha^{-1}}}+\frac{1}{k}\||D|\theta_k\|^2_{H^s_{a,\alpha^{-1}}}+\frac{1}{2}\||D|^\alpha \theta_k\|^2_{H^s_{a,\alpha^{-1}}} \leq C  \|\theta_k\|^4_{H^s_{a,\alpha^{-1}}},$$
Integrating on  $[0, t]$, we obtain
$$\|\theta_k(t)\|^2_{H^s_{a,\alpha^{-1}}}+\frac{1}{k}\int^t_0 \||D|\theta_k\|^2_{H^s_{a,\alpha^{-1}}}+\int^t_0 \||D|^{\alpha}\theta_k\|^2_{H^s_{a,\alpha^{-1}}}
\leq \|\theta^0\|^2_{H^s_{a,\alpha^{-1}}}+ C \int^t_0\|\theta_k(t)\|^4_{H^s_{a,\alpha^{-1}}}$$
Let $t_k >0$ defined by  $t_k = \sup \{t\geq 0 /\sup_{0\leq z\leq t} \|\theta_k(z)\|_{H^s_{a,\alpha^{-1}}}<2  \|\theta^0\|_{H^s_{a,\alpha^{-1}}}\}.$\\
By continuity of $(t\mapsto \|\theta_k(t)\|_{H^s_{a,\alpha^{-1}}})$ we get $0<t_k<T^*_k$.\\
 For all $0\leq t< t_k$, we have
$$\|\theta_k(t)\|^2_{H^s_{a,\alpha^{-1}}}\leq \|\theta^0\|^2_{H^s_{a,\alpha^{-1}}}+ C t\|\theta^0\|^{4}_{H^s_{a,\alpha^{-1}}},$$
 then $$\|\theta_k(t)\|^2_{H^s_{a,\alpha^{-1}}}+\frac{1}{k}\int^t_0 \||D|\theta_k\|^2_{H^s_{a,\alpha^{-1}}}+\int^t_0 \||D|^{\alpha}\theta_k\|^2_{H^s_{a,\alpha^{-1}}}\leq \|\theta^0\|^2_{H^s_{a,\alpha^{-1}}}+ 16 C t \|\theta^0\|^4_{H^s_{a,\alpha^{-1}}}$$
Let  $T>0$ verfies $1+16 T\|\theta^0\|^2_{H^s_{a,\alpha^{-1}}} = 2:$\\
then  $$T= \frac{1}{C\|\theta^0\|^{2}_{H^s_{a,\alpha^{-1}}}}.$$
For $t\in [0, \min (T, t_k))$ one has
$$ \|\theta_k(t)\|^2_{H^s_{a,\alpha^{-1}}}\leq 2 \|\theta^0\|^2_{H^s_{a,\alpha^{-1}}}.$$
Therefore, by continuity of $(t\mapsto \|\theta(t)\|_{H^s_{a,\alpha^{-1}}})$, we get $t_k>T.$\\
 To choose the above $T$ one deduce that  $(\theta_k) $ is bounded in $C([0, T],H^s_{a,\alpha^{-1}} )$ and $(|D|^\alpha \theta_k )$ in $ L^2([0, T], H^{s}_{a,\alpha^{-1}})$. Moreover, we have, for $0\leq t\leq T$
 \begin{align}\label{Mt}
 \|\theta_k(t)\|^2_{H^s_{a,\alpha^{-1}}}+\frac{1}{k}\int^t_0 \||D|\theta_k\|^2_{H^s_{a,\alpha^{-1}}}+\int^t_0 \||D|^{\alpha}\theta_k\|^2_{H^s_{a,\alpha^{-1}}}
\leq \|\theta^0\|^2_{H^s_{a,\alpha^{-1}}}+ 16 C t \|\theta^0\|^4_{H^s_{a,\alpha^{-1}}}=M_t.
\end{align}
which is necessarily in the last step.\\
\noindent \textbf{Step3: Uniqueness and local existence of solution for $(QG)$}. Let $\theta_k$ and $\theta_{k'}$ two solutions for $(QG)_k$, $(QG)_{k'}$  respectively with $k<k'$. Put $\omega_{k,k'}=\theta_{k}-\theta_{k'}:\omega$, we get
\begin{equation}\label{diff}
\partial_{t} \omega-\frac{1}{k} \Delta \omega+\left(\frac{1}{k}-\frac{1}{k^{\prime}}\right) \Delta \theta_{k^{\prime}}+(-\Delta)^{\alpha} \omega+u_{\omega} . \nabla \theta_{k}-u_{\theta_{k} \prime} \nabla \omega=0.
\end{equation}
Taking the product scalar in $H^{s-1}_{{a,\alpha^{-1}}} $ with $\omega$ in (\ref{diff}), we get
$$\frac{1}{2}\frac{d}{dt}\|\omega\|^2_{H^{s-1}_{a, \alpha^{-1}}}+\| |D|^{\alpha}\omega\|^2_{H^{s-1}_{a, \alpha^{-1}}}\leq I_1+I_2+I_3+I_4,$$
with
$$\left\{\begin{array}{lcl}
I_1&=&\frac{1}{k}|\langle \Delta \omega,\omega\rangle_{H^{s-1}_{a, \alpha^{-1}}}|,\\
I_2&=&(\frac{1}{k}-\frac{1}{k'})|\langle \Delta \theta_{k'},\omega\rangle_{H^{s-1}_{a, \alpha^{-1}}}|,\\
I_3&=&| \langle u_{\omega} . \nabla \theta_{k}, \omega \rangle_{H^{s-1}_{a, \alpha^{-1}}}|,\\
I_4&=&| \langle u_{\theta_{k} \prime}. \nabla \omega, \omega\rangle_{H^{s-1}_{a, \alpha^{-1}}}|
\end{array}\right.$$
Then we get
$$I_1\leq\frac{1}{k}\|\nabla\omega\|^2_{H^{s-1}_{a, \alpha^{-1}}}.$$
and
$$I_2\leq(\frac{1}{k}-\frac{1}{k'})\|\nabla \theta_{k'}\|_{H^{s-1}_{a, \alpha^{-1}}}\|\nabla\omega\|_{H^{s-1}_{a, \alpha^{-1}}}.$$
Since $s>1$ we get $H^{s-1}_{a, \alpha^{-1}}$ is algebra, then
$$I_3\leq C_1 \| \omega \|^2_{H^{s-1}_{a, \alpha^{-1}}} \| \nabla \theta_{k} \|_{H^{s-1}_{a, \alpha^{-1}}}$$
Using the lemma \ref{SG} (\ref{eql12}), to estimate $I_4$
\begin{align*}
I_4&\leq C_2 \| \theta_{k'}\|_{H^{s-1}_{a, \alpha^{-1}}}\| |D|^{\alpha} \omega\|_{H^{s-1}_{a, \alpha^{-1}}} \| \omega\|_{H^{s-1}_{a, \alpha^{-1}}}\\
&\leq 2 C_2\| \theta^0\|_{H^{s}_{a, \alpha^{-1}}}.
\end{align*}
By inequality $ab\leq \frac{a^2}{2}+\frac{b^2}{2}$ for $I_4$, we obtain
\begin{align*}
\frac{1}{2}\frac{d}{dt}\|\omega\|^2_{H^{s-1}_{a, \alpha^{-1}}}+\frac{1}{2}\| |D|^{\alpha}\omega\|^2_{H^{s-1}_{a, \alpha^{-1}}}&\leq \frac{1}{k}\|\nabla\omega\|^2_{H^{s-1}_{a, \alpha^{-1}}}+(\frac{1}{k}-\frac{1}{k'})\| \theta_{k'}\|_{H^{s-1}_{a, \alpha^{-1}}}\|\nabla\omega\|_{H^{s-1}_{a, \alpha^{-1}}}+ C\|\omega\|^2_{H^{s-1}_{a, \alpha^{-1}}}\\
&\leq C(\frac{1}{k}+\frac{1}{k'})\| \theta^0\|^2_{H^{s-1}_{a, \alpha^{-1}}}+C \|\omega\|^2_{H^{s-1}_{a, \alpha^{-1}}}.
\end{align*}
Then $$\frac{1}{2}\frac{d}{dt}\|\omega\|^2_{H^{s-1}_{a, \alpha^{-1}}}\leq C(\frac{1}{k}+\frac{1}{k'})\| \theta^0\|^2_{H^{s-1}_{a, \alpha^{-1}}}+C\|\omega\|^2_{H^{s-1}_{a, \alpha^{-1}}}$$
Using the Gronwall's lemma to obtain for $0<t<T$
\begin{align*}
\|\omega(t)\|^2_{H^{s-1}_{a, \alpha^{-1}}}&\leq C t(\frac{1}{k}+\frac{1}{k'})\| \theta^0\|^2_{H^{s-1}_{a, \alpha^{-1}}}\exp(C t)\\
&\leq C T(\frac{1}{k}+\frac{1}{k'})\| \theta^0\|^2_{H^{s-1}_{a, \alpha^{-1}}}\exp(C') \longrightarrow 0 \ if \ k, k'\rightarrow \infty.
\end{align*}
Thus, $\lim_{k\rightarrow\infty}\theta_k $ exists:$\theta \in \mathcal{C}\left([0, T] ; H^{s-1}_{a, \alpha^{-1}}\right)$.\\
Moreover, we have $$\|\omega(t)\|^2_{H^{s-1}_{a, \alpha^{-1}}}+\int^T_0 \| |D|^{\alpha} \omega\|^2_{H^{s-1}_{a, \alpha^{-1}}}\leq C T(\frac{1}{k}+\frac{1}{k'})\| \theta^0\|^2_{H^{s-1}_{a, \alpha^{-1}}}\exp(C'\|\theta^0\|_{H^{s-1}_{a, \alpha^{-1}}}) \longrightarrow 0 \ when \ k, k'\rightarrow \infty.$$
Which implies that $(|D|^{\alpha}\theta_k)$ is a Cauchy sequence in the Hilbert space $L^2([0, T], H^{s-1}_{a, \alpha^{-1}})$, then by  uniqueness of limit, we get $|D|^\alpha \theta \in L^2([0, T], H^{s-1}_{a, \alpha^{-1}})$.\\
 \textbf{The uniqueness of solution }\\
Suppose we have two solutions of $(QG)$ $\theta^1$ and $\theta^2$.\\
Let $\omega =\theta ^1-\theta ^2$ in $\mathcal{C}\left([0, T] ; H^{s-1}_{a, \alpha^{-1}}\right)$ then $\omega(0)=0$
and $$\partial_{t} \omega+(-\Delta)^{\alpha} \omega+\delta . \nabla \theta^{1}-u_{\theta^2}. \nabla \omega=0.$$
Taking the scalar product in $L^2$ in  the last equation with $\omega$ we thus get
\begin{align*}
\frac{1}{2}\frac{d}{dt}\|\omega\|^2_{L^2}+\| |D|^{\alpha}\omega\|^2_{L^2}
&\leq | \langle\omega . \nabla \theta^{1}, \omega \rangle|_{L^2}+| \langle u_{\theta^2}. \nabla \omega, \omega\rangle|_{L^2}\\
&\leq C_1 \| \nabla \theta^{1} \|_{L^2}\|\omega\|^2_{L^2}\\
&\leq C_2  \| \nabla \theta^{1} \|_{H^{s-1}_{a, \alpha^{-1}}}\|\omega\|^2_{L^2}.
\end{align*}
By inequality $ab\leq \frac{a^2}{2}+\frac{b^2}{2}$, we get
$$\frac{d}{dt}\|\omega\|^2_{L^2}+2\| |D|^{\alpha}\omega\|^2_{L^2}\leq  C \|\omega\|^2_{L^2}$$
Since $\omega (0)=0$ Gronwall Lemma gives us the desired result.\\
\textbf{Step4: Continuity in time of the solution.}\\
$\bullet$ We have $\theta_k \rightharpoonup \theta$ uniformly in $H^{s}_{a, \alpha^{-1}}$. Indeed:\\
Since $\theta_k \rightarrow \theta$ in $ H^{s-1}_{a, \alpha^{-1}}$ then $ \theta_k \rightharpoonup \theta$ in $ H^{s-1}_{a, \alpha^{-1}}$, from remark (2.5) it follows that
\begin{align*}
&\sup_{0\leq t\leq T} |\langle\theta_k(t)/f\rangle_{H^{s-1}_{a, \alpha^{-1}}}-\langle\theta(t)/f\rangle_{H^{s-1}_{a, \alpha^{-1}}}|\rightarrow0,\quad \forall f \in H^{s-1}_{a, \alpha^{-1}}(\mathbb{R}^2)\\
&\sup_{0\leq t\leq T} |\langle\theta_k(t)/g\rangle_{H^{s-1}_{a, \alpha^{-1}}}-\langle\theta(t)/g\rangle_{H^{s-1}_{a, \alpha^{-1}}}|\rightarrow0,\quad \forall g \in e^{-a|D|^{\alpha}}\mathcal{S}(\mathbb{R}^2)
\end{align*}
Let $\varphi =\mathcal{F}^{-1}(\frac{1}{1+|\xi|^2}\widehat{g}(\xi))$
\begin{align*}
&\sup_{0\leq t\leq T} |\langle\theta_k(t)/\varphi\rangle_{H^{s}_{a, \alpha^{-1}}}-\langle\theta(t)/\varphi \rangle_{H^{s}_{a, \alpha^{-1}}}|\rightarrow0,\forall \quad \varphi \in e^{-a|D|^{\alpha}}\mathcal{S}(\mathbb{R}^2)\\
&\sup_{0\leq t\leq T} |\langle\theta_k(t)/\varphi\rangle_{H^{s}_{a, \alpha^{-1}}}-\langle\theta(t)/\varphi \rangle_{H^{s}_{a, \alpha^{-1}}}|\rightarrow0,\forall \quad \varphi \in H^{s}_{a, \alpha^{-1}}(\mathbb{R}^2).
\end{align*}
$\bullet$ We have $\theta \in  L^{\infty}\left([0, T] ; H^{s}_{a, \alpha^{-1}} \right)$. Indeed: Let $\varphi \in H^{s}_{a, \alpha^{-1}}$ such that$ \|\varphi\|_{H^{s}_{a, \alpha^{-1}}}\leq 1$\\
$$\langle\theta(t)/\varphi\rangle_{H^{s}_{a, \alpha^{-1}}}=\langle\theta(t)-\theta_k(t)/\varphi\rangle_{H^{s}_{a, \alpha^{-1}}}+\langle\theta_k(t)/\varphi\rangle_{H^{s}_{a, \alpha^{-1}}}$$then
\begin{align*}
|\langle\theta(t)/\varphi\rangle_{H^{s}_{a, \alpha^{-1}}}|&\leq|\langle\theta(t)-\theta_k(t)/\varphi\rangle_{H^{s}_{a, \alpha^{-1}}}|+|\langle\theta_k(t)/\varphi\rangle_{H^{s}_{a, \alpha^{-1}}}|\\
&\leq\sup_{0\leq z\leq T} |\langle\theta(z)-\theta_k(z)/\varphi\rangle_{H^{s}_{a, \alpha^{-1}}}|+\sup_{0\leq z\leq T}|\langle\theta_k(t)/\varphi\rangle_{H^{s}_{a, \alpha^{-1}}}|\\
&\leq\sup_{0\leq z\leq T} |\langle\theta(z)-\theta_k(z)/\varphi\rangle_{H^{s}_{a, \alpha^{-1}}}|+2\|\theta^0\|_{H^{s}_{a, \alpha^{-1}}}
\end{align*}
$k\rightarrow\infty$ implies
$$|\langle\theta(t)/\varphi\rangle_{H^{s}_{a, \alpha^{-1}}}|\leq 2\|\theta^0\|_{H^{s}_{a, \alpha^{-1}}}$$
then $$\sup_{\|\varphi\|_{H^{s}_{a, \alpha^{-1}}}\leq 1}|\langle\theta(t)/\varphi\rangle_{H^{s}_{a, \alpha^{-1}}}|=\|\theta (t)\|_{H^{s}_{a, \alpha^{-1}}}\leq 2\|\theta^0\|_{H^{s}_{a, \alpha^{-1}}}.$$
$\bullet$ $ \theta $ is solution of $(QG)$. In fact\\
$$\partial_t \theta_k\rightharpoonup \partial_t\theta.$$
$$ |D|^{2\alpha}\theta_k\rightharpoonup |D|^{2\alpha}\theta.$$
$$ u_{\theta_k}.\nabla \theta_k\rightharpoonup u_{\theta}.\nabla \theta,$$ in fact
\begin{align*}
\|u_{\theta_k}.\nabla \theta_k- u_{\theta}.\nabla \theta\|_{H^{s-2}_{a,\alpha^{-1} }}&\leq \|(u_{\theta_k}-u_{\theta}) \theta_k- u_{\theta} (\theta_k-\theta)\|_{H^{s-1}_{a,\alpha^{-1} }}\\
&\leq\|(u_{\theta_k}-u_{\theta}) \theta_k\|_{H^{s-1}_{a,\alpha^{-1} }}+\| u_{\theta} (\theta_k-\theta)\|_{H^{s-1}_{a,\alpha^{-1} }}\\
&\leq \|(u_{\theta_k}-u_{\theta})\|_{H^{s-1}_{a,\alpha^{-1} }}\|\theta\|_{H^{s}_{a,\alpha^{-1} }}\rightarrow 0.
\end{align*}
$\bullet \ \theta$ is weakly continuous in $H^{s}_{a, \alpha^{-1}}$. Indeed,
let $\varphi \in H^s_{a, \alpha^{-1}}$ and $0<t\leq t'\leq T$, we have
\begin{align*}
|\langle\theta(t)-\theta(t')/\varphi\rangle_{H^{s}_{a, \alpha^{-1}}}|&\leq |\langle\theta(t)-\theta_k(t)+\theta_k(t)-\theta_k(t')+\theta_k(t')-\theta(t')/\varphi\rangle_{H^{s}_{a, \alpha^{-1}}}|\\
&\leq 2 \sup_{0\leq z\leq T} |\langle\theta(z)-\theta_k(z)/\varphi\rangle_{H^{s}_{a, \alpha^{-1}}}|+|\langle\theta_k(t)-\theta_k(t')/\varphi\rangle_{H^{s}_{a, \alpha^{-1}}}|\rightarrow_{t\mapsto t'}0.
\end{align*}
$\bullet \ \theta $ is right-continuous at $t_0$. Indeed,
Since $\theta$ is weakly continuous in $H^{s}_{a, \alpha^{-1}}$ it suffices to prove
$$ \limsup_{t\rightarrow 0^+}\|\theta (t)\|_{H^{s}_{a, \alpha^{-1}}}\leq \|\theta^0\|_{H^{s}_{a, \alpha^{-1}}}$$
For $t>0$ the fact ($\theta _k (t)\rightarrow \theta (t)$) in $H^{s-1}_{a, \alpha^{-1}}$ implies the existence a subsequence of $(\theta_k)$ such that $$\widehat{\theta }_{\varphi(k)}\rightarrow \widehat{\theta} \ a. e$$
For (\ref{Mt}), we have $$\|\theta_{\varphi_t(k)}(t)\|^2_{H^{s}_{a, \alpha^{-1}}}\leq M_t. $$
Fatou's Lemma gives
$$ \|\theta (t)\|_{H^{s}_{a, \alpha^{-1}}}\leq \liminf\|\theta_{\alpha(k)}(t)\|_{H^{s}_{a, \varphi_t(k)^{-1}}}.$$
Then
 $\|\theta (t)\|^2_{H^{s}_{a, \alpha^{-1}}}\leq M_t$ and
 \begin{align*}
 &\limsup_{t\rightarrow 0^+}\|\theta (t)\|^2_{H^{s}_{a, \alpha^{-1}}}\leq \|\theta^0\|^2_{H^{s}_{a, \alpha^{-1}}}.
\end{align*}
In sum, $\theta $ is right-continuous at 0.\\
To prove the same result at any $t_0 \in [0, T]$, we consider the follow system
\begin{align*}
   (QG)\hspace {2cm}
   \begin{cases}
   \;\partial_t \gamma+(-\Delta)^{\alpha}\gamma +u_{\gamma}.\nabla \gamma &=0\\
   \;\gamma(0, x) &=\theta(t_0, x).
   \end{cases}
   \end{align*}
By the result just proved, $\gamma $ is continuous at $t=0$ by the uniqueness of solution, we obtain the right-continuous at $t_0$ of $\theta $.\\
$\bullet \ \theta $ is left-continuous at $t_0$. Indeed,
For $0<t \leq t_0$ and $1<s_k<s$ such that $s_k\nearrow s$. One has
$$\|\theta(t_0)\|^2_{H^{s_k}_{a, \alpha^{-1}}}+\frac{1}{2}\int ^{t_0}_t \| |D|^{\alpha }\theta\|^2_{H^{s_k}_{a, \alpha^{-1}}}=\|\theta (t)\|^2_{H^{s_k}_{a, \alpha^{-1}}}-\underbrace{\int^{t_0}_t\langle u_{\theta}.\nabla\theta, \theta\rangle_{H^{s_k}_{a,\alpha^{-1}}}}_I$$
But
\begin{align*}
|I|&\leq 2^{s_k+\alpha-1}(s_k+a\alpha)\int ^{t_0}_t \| \theta\|^2_{H^{s_k}_{a, \alpha^{-1}}}\| |D|^{\alpha }\theta\|_{H^{s_k}_{a, \alpha^{-1}}}\\
&\leq2^{s+\alpha-1}(s+a\alpha)\int ^{t_0}_t \| \theta\|^2_{H^{s_k}_{a, \alpha^{-1}}}\| |D|^{\alpha }\theta\|_{H^{s_k}_{a, \alpha^{-1}}}\\
\end{align*}
Then
$$|I|\leq C(t_0-t)\|\theta^0\|^4_{H^{s_k}_{a, \alpha^{-1}}}+\frac{1}{2}\int ^{t_0}_t \| |D|^{\alpha }\theta\|^2_{H^{s_k}_{a, \alpha^{-1}}}$$ Then
\begin{align*}
\|\theta(t)\|^2_{H^{s_k}_{a, \alpha^{-1}}}&\leq \|\theta (t_0)\|^2_{H^{s_k}_{a, \alpha^{-1}}}+C(t_0-t)\|\theta^0\|^4_{H^{s_k}_{a, \alpha^{-1}}}+\int ^{t_0}_t \| |D|^{\alpha }\theta\|^2_{H^{s_k}_{a, \alpha^{-1}}}\\
&\leq \|\theta (t_0)\|^2_{H^{s}_{a, \alpha^{-1}}}+C(t_0-t)\|\theta^0\|^4_{H^{s}_{a, \alpha^{-1}}}+\int ^{t_0}_t \| |D|^{\alpha }\theta\|^2_{H^{s}_{a, \alpha^{-1}}}.
\end{align*}
By Monotone convergence theorem, we get
$$\|\theta(t)\|^2_{H^{s}_{a, \alpha^{-1}}}\leq \|\theta (t_0)\|^2_{H^{s}_{a, \alpha^{-1}}}+C(t_0-t)\|\theta^0\|^4_{H^{s}_{a, \alpha^{-1}}}+\int ^{t_0}_t \| |D|^{\alpha }\theta\|^2_{H^{s}_{a, \alpha^{-1}}}$$
 Since, $|D|^{\alpha }\theta$ is bounded in $L^2([0, T], H^{s}_{a, \alpha^{-1}})$. Then, we pass $\limsup_{t\mapsto t_0} $ we obtain $$\|\theta(t)\|^2_{H^{s}_{a, \alpha^{-1}}}\leq \|\theta(t_0)\|^2_{H^{s}_{a, \alpha^{-1}}}.$$
 Then, by proposition \ref{thhb}, we get the left continuous at $t_0$ which complete the proof of theorem.
\subsection{Exponential type explosion}\label{sub 6.2}
\begin{theorem}\label{Thlimsup}
Let $\theta^0\in H^s_{a,\alpha^{-1}}$, then there exists a unique maximal solution in $\mathcal{C}([0,T^*),H^s_{a,\alpha^{-1}})$, with  $T^*\in(0,\infty]$.\\
Moreover, if $T^*<\infty$ then
\begin{align}\label{Rlimsup}
\limsup_{t\nearrow T^*}\|\theta(t)\|^2_{H^s_{a,\alpha^{-1}}}=\infty.
\end{align}
\end{theorem}
{\bf Proof.}
This proof is inspired from the work of Nader Masmoudi (see \cite{NM}).\\
By the local existence step of the proof of theorem \ref{TSG}, we get
$$\theta\in \mathcal{C}([0,T_1],H^s_{a,\alpha^{-1}}),$$
with $$T_1=\frac{C}{\|\theta^0\|^2_{H^s_{a,\alpha^{-1}}}}.$$
 Consider the following system
\begin{align*}
   (P_1)\hspace {2cm}
   \begin{cases}
   \;\partial_t \gamma+|D|^{2\alpha}\gamma+u_{\gamma}.\nabla \gamma &=0\\
  %% \;u_{\theta} =(u^1_{\theta}, u^2_{\theta}) &=(-\partial_2|D|^{-1}\theta, \partial_1|D|^{-1}\theta)\\
   \;\gamma(0, x) &=\theta(T_1) \in H^s_{a,\alpha^{-1}}.
   \end{cases}
   \end{align*}
Then there exist unique solution $\gamma\in \mathcal{C}([0,T_2],H^s_{a,\alpha^{-1}})$ of $(P_1)$ such that:
$$T_2=\frac{C}{\|\gamma^0\|^2_{H^s_{a,\alpha^{-1}}}}=\frac{C}{|\theta(T_1)\|^2_{H^s_{a,\alpha^{-1}}}}.$$
By the uniqueness of solution we get: $$\gamma(t)=\theta(T_1+t); \quad \forall 0\leq t\leq  T_2$$
then $\theta\in \mathcal{C} ([0,T_1+T_2],H^s_{a,\alpha^{-1}})$ such that $T_1+T_2<T^*$.
Then we can obtain $T_1,T_2, T_3,..,T_n$ such that
$$T_k= \frac{C}{\|\theta(T_1+..+T_{k-1})\|^2_{H^s_{a,\alpha^{-1}}}}; \quad \forall 1\leq k\leq n.$$
 Consider the following system
\begin{align*}
   (P_k)\hspace {2cm}
   \begin{cases}
   \;\partial_t \delta+|D|^{2\alpha}\delta +u_{\delta}.\nabla \delta &=0\\
  %% \;u_{\theta} =(u^1_{\theta}, u^2_{\theta}) &=(-\partial_2|D|^{-1}\theta, \partial_1|D|^{-1}\theta)\\
   \;\delta(0, x) &=\theta( \sum_{k=1}^n T_k)\in H^s_{a,\alpha^{-1}}.
   \end{cases}
   \end{align*}
Then there is a unique solution of $(P_k)$, $\delta \in \mathcal{C}([0,T_{n+1}],H^s_{a,\alpha^{-1}})$.\\
with
$$T_{n+1}=\frac{C}{\|\theta(\sum_{k=1}^n T_k)\|^2_{H^s_{a,\alpha^{-1}}}}.$$
Also, by uniqueness, we get $$\delta(t)= \theta(\sum_{i=1}^n T_i+t); \quad \forall 0\leq t\leq T_{n+1}.$$
Then
  $$T_1+..+T_{n+1}<T^*,$$
 which gives that $\sum_{k=1}^\infty T_k \leq T^*.$\\
  Then $$T_k\underset{k\rightarrow +\infty}{\longrightarrow} 0.$$
Thus
 $$\frac{C}{\|\theta(\sum_{i=1}^k T_i)\|^2_{H^s_{a,\alpha^{-1}}}}\underset{k\rightarrow \infty}{\longrightarrow} 0,$$
which implies that
 $$\|\theta(\sum_{i=1}^k T_i)\|^2_{H^s_{a,\alpha^{-1}}}\underset{k\rightarrow \infty}{\longrightarrow} \infty.$$
 Then $\sum_{i=1}^\infty T_i=T^*$ and $\limsup_{t\nearrow T^*}\|\theta(t)\|^2_{H^s_{a,\alpha^{-1}}}=+\infty.$
\begin{prop}
Let $s>2 $,  $ a>0$ and $0<\alpha<1/2$. If $ \theta$ is a maximal solution of $(QG)$ in $\mathcal{C}([0, T^*), H^s_{a,\alpha^{-1}}(\mathbb{R}^2)) $ with $T^*$ is finite, then
\begin{align}\label{4.3.1}
\int^{T^*}_t \|\mathcal{F}(e^{a\alpha|D|^{\alpha}}\nabla\theta(\tau))\|^2_{L^1}d\tau=\infty, \quad \ \forall 0 \leq t < T^*.
\end{align}
\begin{align}\label{4.3.2}
\frac{C(a, s)}{T^*-t}\leq \|\theta(t)\|^2_{H^s_{a,\alpha^{-1}}}, \quad \ \forall t\in[0,T^*).
\end{align}
\end{prop}
{\bf Proof.}
Clearly, by theorem \ref{Thlimsup}, $\underset{t\rightarrow T^*}\lim\|\theta(t)\|_{H^s_{a,\alpha^{-1}}}=\infty.$\\
\underline{Proof of \ref{4.3.1}}.
Taking the scalar product with $\theta$ in $H^s_{a,\alpha^{-1}}(\mathbb{R}^2)$ and using (\ref{eql13}), we obtain
\begin{align*}
\frac{1}{2}\frac{d}{dt}\|\theta(t)\|^2_{H^s_{a,\alpha^{-1}}}+\||D|^{\alpha} \theta(t)\|^2_{H^s_{a,\alpha^{-1}}}&\leq|\langle u_{\theta}.\nabla\theta, \theta\rangle_{H^s_{a,\alpha^{-1}}}|\\
&\leq C \|\mathcal{F}(e^{a\alpha|D|^{\alpha}}\nabla\theta)(t)\|_{L^1}\||D|^\alpha\theta\|_{H^s_{a,\alpha^{-1}}}\|\theta\|_{H^s_{a,\alpha^{-1}}},
\end{align*}
with $C= C(a,s,\alpha)=2^{s+\alpha-1}(s+a\alpha)$
 \noindent From inequality $xy \leq \frac{x^2}{2}+\frac{y^2}{2}$ and integrating on $[t, T] \subset [0, T^*)$ we get
$$\|\theta(T)\|^2_{H^s_{a,\alpha^{-1}}}+\int^T_t\||D|^{\alpha} \theta(t)\|^2_{H^s_{a,\alpha^{-1}}}\leq \|\theta(t)\|^2_{H^s_{a,\alpha^{-1}}}+\int^T_t \|\mathcal{F}(e^{a\alpha|D|^{\alpha}}\nabla\theta)(\tau)\|^2_{L^1}\|\theta\|^2_{H^s_{a,\alpha^{-1}}}d\tau$$
with $C(z,s) $ is increasing with respect to $z$ then $C(a\alpha,s)\leq C(a,s)$.\\
Gronwall's Lemma implies
$$\|\theta(T)\|^2_{H^s_{a,\alpha^{-1}}}+\int^T_t\||D|^{\alpha} \theta(\tau)\|^2_{H^s_{a,\alpha^{-1}}}\leq C \|\theta(t)\|^2_{H^s_{a,\alpha^{-1}}} \exp( C \int^T_t\|\mathcal{F}(e^{a\alpha|D|^{\alpha}}\nabla\theta)(\tau)\|^2_{L^1})$$
The fact that
the fact that $\limsup_{T\nearrow T^*}\|\theta(t)\|_{H^s_{a,\alpha^{-1}}}=\infty $ implies
$$\int^{T^*}_t{\|\mathcal{F}(e^{a\alpha|D|^\alpha}\nabla\theta)(t)\|^2_{L^1}d\tau}=\infty,\quad \forall \ 0 \leq t< T^*.$$
and $$\limsup_{T\nearrow T^*}\|\mathcal{F}(e^{a\alpha|D|^\alpha}\nabla\theta)(t)\|^2_{L^1}=\infty.$$
\underline{Proof of \ref{4.3.2}}.
 Due to (2.7) and using inequality $xy \leq \frac{x^2}{2}+\frac{y^2}{2}$,\\
with $ x= \||D|^{\alpha} \theta(t)\|_{H^s_{a,\alpha^{-1}}} $ and $ y = C \|\theta(t)\|^2_{H^s_{a,\alpha^{-1}}}$, we get
 $$\frac{1}{2}\frac{d}{dt}\|\theta(t)\|^2_{H^s_{a,\alpha^{-1}}}+\||D|^{\alpha} \theta(t)\|^2_{H^s_{a,\alpha^{-1}}}\leq \||D|^{\alpha} \theta(t)\|^2_{H^s_{a,\alpha^{-1}}}\|\theta(t)\|^2_{H^s_{a,\alpha^{-1}}}$$
 Integrating  on $[t_0, t] \subset [0, T^*)$ and applying Gronwall's Lemma we get
     $$\|\theta(t)\|^2_{H^s_{a,\alpha^{-1}}}+\int^t_{t_0}\||D|^{\alpha} \theta(\tau)\|^2_{H^s_{a,\alpha^{-1}}}\leq \|\theta(t_0)\|^2_{H^s_{a,\alpha^{-1}}} \exp( C \int^t_{t_0}\|\theta(\tau)\|^2_{H^s_{a,\alpha^{-1}}}).$$
Then
$$\|\theta(t)\|^2_{H^s_{a,\alpha^{-1}}}\exp(-C \int^t_{t_0}\|\theta(\tau)\|^2_{H^s_{a,\alpha^{-1}}})\leq  \|\theta(t_0)\|^2_{H^s_{a,\alpha^{-1}}}.$$
Integrating over $[t_0, T] \subset [0, T^*)$ we obtain
$$1-\exp(-C \int^T_{t_0}\|\theta(\tau)\|^2_{H^s_{a,\alpha^{-1}}})\leq C\|\theta(t_0)\|^2_{H^s_{a,\alpha^{-1}}}(T-t_0).$$
 Taking $T\rightarrow T^*$, we get
$$1-\exp(-C \int^{T^*}_{t_0}\|\theta(\tau)\|^2_{H^s_{a,\alpha^{-1}}})\leq C\|\theta(t_0)\|^2_{H^s_{a,\alpha^{-1}}}(T^*-t_0).$$
But, from Cauchy-Schwartz inequality, we have
$$\|\mathcal{F}(e^{a\alpha|D|^\alpha}\nabla\theta)(t)\|^2_{L^1}\leq C(a, s, \alpha) \|\theta(t)\|^2_{H^s_{a,\alpha^{-1}}}$$
with $$C^2=C^2(a, s, \alpha)=\int e^{-2a(1-\alpha)|\xi|^{\alpha}}|\xi|^2(1+|\xi|^2)^{-s}d\xi <\infty$$
According this with $$\int^{T^*}_t{\|\mathcal{F}(e^{a\alpha|D|^\alpha}\nabla\theta)(\tau)\|_{L^1}d\tau}=\infty,\quad \forall \ 0 \leq t< T^*.$$ we obtain the desired result and Proof of proposition is finished.\\
{\bf Proof of (\ref{TSG.1}).} This proof is done in two steps:\\
$\bullet$ \textbf{Step1} we prove $$\frac{C(a, s)}{T^*-t}\leq \|\theta(t)\|^2_{\dot{H}^s}$$
Since $0<a\sqrt{\alpha}\leq a$ we have the following embedding $H^s_{a,\alpha^{-1}}\hookrightarrow H^s_{a\sqrt{\alpha},\alpha^{-1}}$ which implies that $ \theta\in \mathcal{C}([0,T^*), H^s_{a\sqrt{\alpha},\alpha^{-1}})$ and if $T^*_{a\sqrt{\alpha},\alpha^{-1}}$ the maximal time of existence some solution in $ H^s_{a\sqrt{\alpha},\alpha^{-1}}$ we get
\begin{align}\label{*}
T^* \leq T^*_{a\alpha,\alpha^{-1}}.
\end{align}
The fact that $$\limsup_{T\nearrow T^*}\|\mathcal{F}(e^{a\alpha|D|^\alpha}\nabla\theta)(t)\|^2_{L^1}=\infty$$ implies
\begin{align}\label{**}
 \limsup_ {t\rightarrow T^*} \|\theta(t)\|_{H^s_{a\alpha,\alpha^{-1}}}=\infty
 \end{align}
      Indeed, From Cauchy-Schwartz, we have
      \begin{align*}
      \|\mathcal{F}(e^{a\alpha|D|^{\alpha}}\nabla\theta(\tau))\|_{L^1}&=\int_{\xi} e^{-a(\sqrt{\alpha}-\alpha)|\xi|^{\alpha}}|\xi|(1+|\xi|^2)^{-\frac{s}{2}}(1+|\xi|^2)^{\frac{s}{2}}e^{\sqrt{\alpha}|\xi|^{\alpha}}|\widehat{\theta}(\xi)|d\xi\\
      &\leq C  \|\theta(t)\|_{H^s_{a\sqrt{\alpha},\alpha^{-1}}}
      \end{align*}
Combining (\ref{*}) and (\ref{**}) we obtain $T^* =T^*_{a\sqrt{\alpha},\alpha^{-1}}.$ Therefore,
$$\frac{C(a, s)}{T^*-t}\leq \|\theta(t)\|^2_{H^s_{a\sqrt{\alpha},\alpha^{-1}}},  \quad \ \forall t\in[0,T^*).$$
The same principle gives us  $$T^* =T^*_{a\alpha^{\frac{1}{2}},\alpha^{-1}}=...= T^*_{a\alpha^{\frac{n}{2}}, \alpha^{-1}}.$$
By induction, we can conclude that, for any $n\in \mathbb{N}. $
$$\frac{C(a, s)}{T^*-t}\leq \|\theta(t)\|^2_{H^s_{a\alpha^{\frac{n}{2}},\alpha^{-1}}}$$
 Thanks to the convergence dominated theorem, we get
 $$\frac{C(a, s)}{T^*-t}\leq \|\theta(t)\|^2_{H^s}$$
 As  ${\lim}\underset{t\rightarrow T^*} \|\theta(t)\|_{H^s}=\infty$ and $\|\theta(t)\|_{L^2}\leq\|\theta^0\|_{L^2},  \quad \ \forall t\in[0,T^*).$ Then, there is a time $T_0$ such that $$ \|\theta(t)\|^2_{\dot{H}^s}\geq \|\theta^0\|^2_{L^2}; \forall \ T_0 \leq t< T^*).$$
Then
\begin{align}\label{Hs}
 \frac{C/2}{T^*-t}\leq \|\theta(t)\|^2_{\dot{H}^s}; \quad \forall\  T_0 \leq t< T^*.
 \end{align}
 But, we know that $$\|\theta(t)\|^2_{\dot{H}^s_{a,\alpha^{-1}}} = \sum_{k\geq 0}\frac{(2a)^k}{k!}\|\theta(t)\|^2_{\dot{H}^{s+\frac{k\alpha}{2}}}$$
 $\bullet$ \textbf{Step2}, we prove the exponential type explosion.\\
  Using the energy estimate and interpolation for $0< s< s+\frac{k\alpha}{2}$  we obtain
  \begin{align*}
 &\|\theta(t)\|_{\dot{H}^s}\leq (2\pi)\|\theta(t)\|^{\frac{k\alpha}{2s+k\alpha}}_{L^2}\|\theta(t)\|^{\frac{2s}{2s+k\alpha}}_{\dot{H}^{s+\frac{k\alpha}{2}}}\\
 &\|\theta(t)\|_{\dot{H}^s}\leq(2\pi) \|\theta^0\|^{\frac{k\alpha}{2s+k\alpha}}_{L^2}\|\theta(t)\|^{\frac{2s}{2s+k\alpha}}_{\dot{H}^{s+\frac{k\alpha}{2}}}.
 \end{align*}
  %%%%%%%%%%%%%%%%%%%%%%%%%%%%%%%%%%%%%%%%%%%%%%%%%%%%%%
  From (\ref{Hs}), we get
   $$\frac{C/2}{T^*-t}\leq(2\pi)\|\theta^0\|^{\frac{k\alpha}{2s+k\alpha}}_{L^2}\|\theta(t)\|^{\frac{2s}{2s+k\alpha}}_{\dot{H}^{s+\frac{k\alpha}{2}}},$$
  and
   $$(\frac{C/2}{T^*-t})^{\frac{2s+k\alpha}{s}}\leq(2\pi)\|\theta^0\|^{\frac{k\alpha}{s}}_{L^2}\|\theta(t)\|^2_{\dot{H}^{s+\frac{k\alpha}{2}}},$$
   which implies
   $$\frac{C_1}{(T^*-t)^2}. (\frac{C_2}{T^*-t})^{\frac{k\alpha}{s}}\leq \|\theta(t)\|^2_{\dot{H}^{s+\frac{k\alpha}{2}}}$$
   With
   \begin{align*}
   \begin{cases}
   &C_1=(\frac{C}{2})^{2s}\\
   &C_2=\frac{C}{2(2\pi)\|\theta^0\|_{L^2}}
   \end{cases}
   \end{align*}
Multiplying both sides by  $\frac{(2a)^k}{k!}$ we thus get
     $$\frac{C_1}{(T^*-t)^2}. \frac{1}{k!}(\frac{2aC_2}{(T^*-t)})^{\frac{k\alpha}{s}}\leq \frac{(2a)^k}{k!}\int_{\mathbb{R}^2}|\xi|^{2s}|\xi|^{k\alpha}|\widehat{\theta}(\xi)|^2d\xi$$
  Summing over $\{k\geq 0\} $ to obtain $$\frac{C_1}{(T^*-t)^2} \exp(\frac{2aC_2}{(T^*-t)^{\frac{\alpha}{s}}})\leq \|\theta(t)\|^2_{\dot{H}^s_{a,\alpha^{-1}}}.$$
  Which complete the proof.
  \section{Global Solution}
In this section, we give the proof of theorem \ref{TGS}. Let $\theta\in\mathcal{C}((0,T^*),H^s_{a,\alpha^{-1}})$ be a global solution of $(QG)$.\\
$\bullet$ We start by proving that if $T^*<\infty$, then $$\int^{T^*}_0\||D|^{\alpha}\theta\|^2_{H^s_{a,\alpha^{-1}}}d\tau=\infty.$$
%%%%%%%%%%%%%%%%%%%%%%%%%%%%%%%%%%%%%%%%%%%%%%%%%%%%%%%%%%%%%%%%%%%%%%%%%%%%%%%%%%%%%%%%%%%%%%%%%%%%
Assume that there exists $M>0$ such that $\int^{T^*}_0\||D|^{\alpha}\theta\|^2_{H^s_{a,\alpha^{-1}}}d\tau\leq M$\\
For a select time $T_0 \in (0, T^*)$ such that $\int^{T^*}_{T_0} \| |D|^{\alpha}\theta\|^2_{H^{2-2\alpha}}d\tau <\frac{1}{2}.$\\
Then, for all $t\in [T_0,T^*)$
\begin{align*}
\|\theta(t)\|^2_{H^s_{a,\alpha^{-1}}}+ \int^t_{T_0}\|\theta\|^2_{H^s_{a,\alpha^{-1}}}d\tau&\leq \|\theta(T_0)\|^2_{H^s_{a,\alpha^{-1}}}+C\int^t_{T_0}\||D|^\alpha\theta(\tau)\|^2_{H^s_{a,\alpha^{-1}}}\||D|^\alpha\theta\|_{H^s_{a,\alpha^{-1}}}d\tau\\
&\leq\|\theta(T_0)\|^2_{H^s_{a,\alpha^{-1}}}+C\sup_{z\in[T_0,t]}\|\theta(z)\|^2_{H^s_{a,\alpha^{-1}}}\int^t_{T_0}\||D|^\alpha\theta\|_{H^s_{a,\alpha^{-1}}}d\tau\\
&\leq\|\theta(T_0)\|^2_{H^s_{a,\alpha^{-1}}}+C\sup_{z\in[T_0,t]}\|\theta(z)\|^2_{H^s_{a,\alpha^{-1}}}(\sqrt{t-T_0}\int^t_{T_0}\||D|^\alpha\theta\|^2_{H^s_{a,\alpha^{-1}}}d\tau)^{1/2}\\
&\leq\|\theta(T_0)\|^2_{H^s_{a,\alpha^{-1}}}+(\sup_{z\in[T_0,t]}\|\theta(z)\|^2_{H^s_{a,\alpha^{-1}}}C\sqrt{T^*-T_0}\sqrt{M}).
\end{align*}
We can suppose that $$C\sqrt{T^*-T_0}\sqrt{M}<\frac{1}{2},$$
Thus, we get
 $$\sup_{z\in[T_0,T^*]}\|\theta(t)\|^2_{\dot{H}^{2-2\alpha}}\leq 2\|\theta(T_0)\|^2_{H^s_{a,\alpha^{-1}}}.$$
Put $$ M_0= \max(\max_{z\in[0,T_0]}\|\theta(z)\|_{H^s_{a,\alpha^{-1}}}; \sqrt{2}\|\theta(T_0)\|_{H^s_{a,\alpha^{-1}}}),$$
we get $$\|\theta(t)\|_{H^s_{a,\alpha^{-1}}}\leq M_0$$
which contradicts with (\ref{Rlimsup}).\\
%%%%%%%%%%%%%%%%%%%%%%%%%%%%%%%%%%%%%%%%%%%%%%%%%%%%%%%%%%%%%%%%%%%%%%%%%%%%%%%%%%%%%%%%%
$\bullet$ Now, we show that if $\theta^0 \in H^s_{a,\alpha^{-1}}(\mathbb{R}^2)$ such that $ \|\theta^0\|^2_{H^s_{a,\alpha^{-1}}}<\frac{1}{4C}$, we get a global solution in $\mathcal{C}_b(\mathbb{R}^+,H^s_{a,\alpha^{-1}}(\mathbb{R}^2))$.\\
Combining (\ref{eql15}) and the convex inequality $ab\leq \frac{a^2}{2}+\frac{b^2}{2}$ we get
 $$\|\theta(t)\|^2_{H^s_{a,\alpha^{-1}}}+ \int^t_0\||D|^{\alpha}\theta\|^2_{H^s_{a,\alpha^{-1}}}d\tau\leq \|\theta^0\|^2_{H^s_{a,\alpha^{-1}}}+C\int^t_0\|\theta(\tau)\|^2_{H^s_{a,\alpha^{-1}}}\||D|^{\alpha}\theta(\tau)\|^2_{H^s_{a,\alpha^{-1}}}d\tau$$
Let
 $$T= \sup_{0\leq t\leq T^*}\{\sup_{0\leq z\leq t}\|\theta(z)\|_{H^s_{a,\alpha^{-1}}}< 2\|\theta^0\|_{H^s_{a,\alpha^{-1}}} \}$$
 By continuity of $(t\rightarrow \|\theta(t)\|_{H^s_{a,\alpha^{-1}}})$ we get $T>0$ and for all $ 0\leq t<T$, we have
$$\|\theta(t)\|^2_{H^s_{a,\alpha^{-1}}}+ \int^t_0\||D|^{\alpha}\theta\|^2_{H^s_{a,\alpha^{-1}}}d\tau\leq \|\theta^0\|^2_{H^s_{a,\alpha^{-1}}}.$$
This implies $$T=T^*,$$
and $$\int^{T^*}_0\||D|^{\alpha}\theta\|^2_{H^s_{a,\alpha^{-1}}}d\tau<\infty.$$
 Hence $T^*= \infty$ and
$$\|\theta(t)\|^2_{H^s_{a,\alpha^{-1}}}+ \int^t_0\||D|^{\alpha}\theta\|^2_{H^s_{a,\alpha^{-1}}}d\tau\leq \|\theta^0\|^2_{H^s_{a,\alpha^{-1}}}, \quad \forall \ t\geq 0.$$
%%%%%%%%%%%%%%%%%%%%%%%%%%%%%%%%%%%%%%%%%%%%%%%%%%%
 \subsection{Long time decay of global solution }\label{sub 6.3}
 In this section, we prove if $\theta \in \mathcal{C}(\mathbb{R}^+, H^s_{a,\alpha^{-1}}(\mathbb{R}^2))$ is a global solution of $(QG)$ then $\|\theta(t)\|^2_{H^s_{a,\alpha^{-1}}}$ decays to zero as time goes to infinity.
As a first step for proving the theorem \ref{Longtime}, let us show the following proposition:
\begin{prop}\label{borne}
If $\theta$ is a solution of $(QG)$ such that $\theta \in \mathcal{C}_b(\mathbb{R}^+, H^s_{a,\alpha^{-1}}(\mathbb{R}^2))$ and $|D|^\alpha \theta \in L^2(\mathbb{R}^+, H^s_{a,\alpha^{-1}}(\mathbb{R}^2))$, then
$$\lim_{t\rightarrow\infty}\|\theta(t)\|^2_{H^s_{a,\alpha^{-1}}}=0.$$                                                                                                \end{prop}
{\bf Proof.}
We recall the following $L^2$ energy estimate
\begin{equation}\label{eq111}
\|\theta\|_{L^2}^2+2\int_0^t\||D|^\alpha\theta\|_{L^2}^2\leq \|\theta^0\|_{L^2}^2.
\end{equation}
As $0\leq \alpha\leq s$, then
$$|\xi|^s\leq \left\{\begin{array}{l}
|\xi|^\alpha\;{\rm if}\;|\xi|<1\\
|\xi|^{s+\alpha}\;{\rm if}\;|\xi|>1
\end{array}\right.\leq (1+|\xi|^2)^{\frac{s}{2}}|\xi|^\alpha.$$
By the above inequality we obtain
\begin{align*}
\|\theta(t)\|^2_{\dot{H}^s_{a,\alpha^{-1}}}&= \int_{\R^2}|\xi|^{2s}e^{2a|\xi|^\alpha}|\widehat{\theta}(\xi) |^2d\xi\\
&\leq \int_{\{|\xi|<1\}}|\xi|^{2s}e^{2a|\xi|^\alpha}|\widehat{\theta}(\xi) |^2d\xi + \int_{\{|\xi|>1\}}|\xi|^{2s}e^{2a|\xi|^\alpha}|\widehat{\theta} (\xi)|^2d\xi \\
& \leq e^{2a}\|\theta(t)\|^2_{\dot{H}^\alpha} + \||D|^\alpha\theta(t)\|^2_{H^s_{a,\alpha^{-1}}}
\end{align*}
Combining (\ref{eq111}) and the fact $|D|^\alpha \theta \in L^2(\mathbb{R}^+, H^s_{a,\alpha^{-1}}(\mathbb{R}^2))$, we obtain $\theta\in L^2(\mathbb{R}^+, \dot{H}^s_{a,\alpha^{-1}}).$
% and $\theta\mapsto \|\theta(t)\|_{\dot{H}^s_{a,\alpha^{-1}}}$.\\
Let $ \varepsilon>0$ and $$\mathcal{A}_{\varepsilon}=\{t\geq0; \|\theta(t)\|_{H^s_{a,\alpha^{-1}}}>\varepsilon\}.$$
 We have
$$ M_0:e^{2a}\|\theta^0\|^2_{L^2}+\|\theta^0\|^2_{L^2} \geq\int_{\R^3}\|\theta(t)\|^2_{H^s_{a,\alpha^{-1}}}\geq \int_{\mathcal{A}_{\varepsilon}}\|\theta(t)\|^2_{H^s_{a,\alpha^{-1}}}\geq \lambda_1(\mathcal{A}_{\varepsilon})\varepsilon^2$$
Let $T_{\varepsilon}= \varepsilon^{-2} M_0$ then
$\lambda(\mathcal{A}_{\varepsilon})\leq T_{\varepsilon}$ and there is a $ t_0 \in [0,T_{\varepsilon}+2] \setminus \mathcal{A}_{\varepsilon}$
which implies
\begin{align}\label{smalldata}
 \|\theta(t_0)\|^2_{\dot{H}^s_{a,\alpha^{-1}}} < \varepsilon.
 \end{align}
Let us consider the following equation
\begin{align*}
   (QG)\hspace {2cm}
   \begin{cases}
   \;\partial_t \gamma+|D|^{2\alpha}\gamma +u_{\gamma}.\nabla \gamma &=0\\
   %\;u_{\theta} =(u^1_{\theta}, u^2_{\theta}) &=(-\partial_2|D|^{-1}\theta, \partial_1|D|^{-1}\theta)\\
   \;\gamma(0) =\theta(t_0).
   \end{cases}
   \end{align*}
The existence and uniqueness of a solution to the quasi geostrophic equation gives for all $t>0, \ \gamma(t)=\theta(t_0+t)$ then
$$\|\theta(t_0+t)\|^2_{H^s_{a,\alpha^{-1}}}+ \int^t_0\||D|^{\alpha}\theta(t_0+\tau)\|^2_{H^s_{a,\alpha^{-1}}}d\tau\leq \|\theta(t_0)\|^2_{H^s_{a,\alpha^{-1}}}<\varepsilon^2, \quad \forall t>0,$$
which completes the proof of proposition.\\
{\bf Proof theorem \ref{Longtime}.}
By the embedding  $ H^s_{a,\alpha^{-1}}\hookrightarrow H^s_{a\alpha ,\alpha^{-1}}$ we get
    $$\theta \in \mathcal{C}(\mathbb{R}^+, H^s_{a\alpha ,\alpha^{-1}}(\mathbb{R}^2))$$
    Similarly, we get for $k\in \mathbb{N}$
    $$\theta \in \mathcal{C}(\mathbb{R}^+, H^s_{a\alpha^k ,\alpha^{-1}}(\mathbb{R}^2))$$
On the other hand, (see the appendix), it is shown that
      $$\lim_{t\rightarrow\infty}\|\theta(t)\|_{H^s}=0,$$ which implies for any small enough $ \varepsilon>0$ there exists a $t_0>0$ such that
      $$ \|\theta(t)\|_{H^s} < \varepsilon, \quad t\geq t_0.$$
    But the Dominated Convergence Theorem yields $ \lim _{k\rightarrow\infty}  \|\theta(t_0)\|_{H^s_{a\alpha^k ,\alpha^{-1}}}=\|\theta(t_0)\|_{H^s}$.
    Hence, there exists $k_0>0$ and for all $ k\geq k_0$ we have
    $$\|\theta(t_0)\|_{H^s_{a\alpha^k ,\alpha^{-1}}}< \varepsilon.$$
    Using proposition \ref{borne}, we get
    $$\|\theta(t_0+t)\|^2_{H^s_{a\alpha^{k_0},\alpha^{-1}}}\leq \varepsilon, \quad \forall t\geq0,$$
    and $$\|\theta(t_0+t)\|^2_{H^s_{a\alpha^{k_0},\alpha^{-1}}}+ \int^t_0\||D|^{\alpha}\theta(t_0+\tau)\|^2_{H^s_{a\alpha^{k_0},\alpha^{-1}}}d\tau\leq \|\theta(t_0)\|^2_{H^s_{a\alpha^{k_0},\alpha^{-1}}}<\varepsilon^2, \quad \forall t>0,$$
   It follows from (\ref{eql13}) of Lemma \ref{SG} that
\begin{align*}
|\langle u_{\theta}.\nabla\theta, \theta \rangle_{H^s_{a,\alpha^{-1}}}|&\leq C\|e^{a\alpha|D|^\alpha}\theta\|_{X^1}\||D|^\alpha\theta\|_{H^s_{a,\alpha^{-1}}}\|\theta\|_{H^s_{a,\alpha^{-1}}}\\
&\leq C C_1\||D|^\alpha\theta\|_{H^s_{a\alpha,\alpha^{-1}}}\||D|^\alpha\theta\|_{H^s_{a,\alpha^{-1}}}\|\theta\|_{H^s_{a,\alpha^{-1}}}
\end{align*}
with
\begin{align*}
C_1=\int\frac{|\xi|^2}{(1+|\xi|^2)^s|\xi|^{2\alpha}}&=\pi\int^\infty_0\frac{r^3}{(1+r^2)^s r^{\alpha}}dr\\
&<\infty \ ( 2s-3+\alpha>1)
\end{align*}
and $CC_1$ is independent of $a$.\\
 In order to bound $\theta(t)$ in $ H^s_{a,\alpha^{-1}}(\mathbb{R}^2)$, we take the scalar product with $\theta$
 in $ H^s_{a\alpha^{k_0-1} ,\alpha^{-1}}$ to get
\begin{align*}
   \frac{1}{2}\frac{d}{dt} \|\theta(t)\|^2_{H^s_{a\alpha^{k_0-1} ,\alpha^{-1}}}+ \||D|^{\alpha}\theta\|^2_{H^s_{a\alpha^{k_0-1} ,\alpha^{-1}}}&\leq C \||D|^\alpha\theta\|_{H^s_{a\alpha^{k_0}, \alpha^{-1}}} \||D|^{\alpha}\theta\|_{H^s_{a\alpha^{k_0-1} ,\alpha^{-1}}}\|\theta\|_{H^s_{a\alpha^{k_0-1},\alpha^{-1}}}\\
   &\leq\frac{1}{2}\||D|^{\alpha}\theta\|^2_{H^s_{a\alpha^{k_0-1} ,\alpha^{-1}}}
   +\frac{C}{2}\|\theta\|^2_{H^s_{a\alpha^{k_0-1},\alpha^{-1}}} \||D|^\alpha\theta\|^2_{H^s_{a\alpha^{k_0}, \alpha^{-1}}}
    \end{align*}
Integrating over $[t_0,t]$ we obtain
\begin{align*}
    \|\theta(t)\|^2_{H^s_{a\alpha^{k_0-1} ,\alpha^{-1}}}+ \int^t_{t_0}\||D|^{\alpha}\theta\|^2_{H^s_{a\alpha^{k_0-1} ,\alpha^{-1}}}d\tau&\leq \|\theta(t_0)\|^2_{H^s_{a\alpha^{k_0-1} ,\alpha^{-1}}}\\
    &+C\int^t_{t_0} \||D|^{\alpha}\theta(\tau)\|^2_{H^s_{a\alpha^{k_0}, \alpha^{-1}}} \|\theta(\tau)\|^2_{H^s_{a\alpha^{k_0-1} ,\alpha^{-1}}}d\tau
    \end{align*}
    By Gronwall's Lemma we get
    \begin{align*}
     \|\theta(t)\|^2_{H^s_{a\alpha^{k_0-1} ,\alpha^{-1}}}+ \int^t_{t_0}\||D|^{\alpha}\theta\|^2_{H^s_{a\alpha^{k_0-1} ,\alpha^{-1}}}d\tau&\leq
    \|\theta(t_0)\|^2_{H^s_{a\alpha^{k_0-1} ,\alpha^{-1}}}\exp\big(C\int^t_{t_0} \||D|^{\alpha}\theta(\tau)\|^2_{H^s_{a\alpha^{k_0}, \alpha^{-1}}}\big)\\
    &\leq \|\theta(t_0)\|^2_{H^s_{a\alpha^{k_0-1} ,\alpha^{-1}}}\exp\big(C\int^\infty_{t_0} \||D|^{\alpha}\theta(\tau)\|^2_{H^s_{a\alpha^{k_0}, \alpha^{-1}}}\big)\\
  &\leq \|\theta(t_0)\|^2_{H^s_{a\alpha^{k_0-1} ,\alpha^{-1}}}\exp\big(C\int^\infty_{0} \||D|^{\alpha}\theta(\tau)\|^2_{H^s_{a\alpha^{k_0}, \alpha^{-1}}}\big)\\
  &\leq C_{k_0}.
  \end{align*}
  Then $$\|\theta(t)\|^2_{H^s_{a\alpha^{k_0-1} ,\alpha^{-1}}}=M_{k_0}, \forall t\geq0,$$
  with $$M_{k_0}=\max(\|\theta\|_{L^\infty([0,t_0], H^s_{a\alpha^{k_0-1} ,\alpha^{-1}}), C_{k_0}}).$$
  Similarly, by finite and decreasing induction, we can deduce that $\theta(t)$ is bounded in $H^s_{a,\alpha^{-1}}(\mathbb{R}^2), \quad t\geq 0$.\\
     Using the Lemma \ref{borne} to complete the proof of theorem.
   \section{Appendix}
   %%%%%%%%%%%%%%%%%%%%%%%%%%%%%%%%
   In this section, we prove that $$\lim_{t\rightarrow\infty}\|\theta(t)\|_{H^s}=0,$$ if $\theta\in \mathcal{C} (\R^+, H^s_{a,\alpha^{-1}}(\mathbb{R}^2))$ is a global solution of $(QG)$.\\
   In \cite{JC}, we have
   \begin{align*}
   \begin{cases}
   &\lim_{t\rightarrow\infty}\|\theta(t)\|_{\dot{H}^{2-2\alpha}}=0,\\
   &\int^\infty_0\||D|^\alpha\theta\|^2_{\dot{H}^{2-2\alpha}}<\infty.
   \end{cases}
   \end{align*}
   Using Lemma \ref{lemjc}, we get for all $t>t_0$
   $$
\|\theta(t)\|^2_{H^s}+2 \int^t_{t_0}\||D|^{\alpha}\theta\|^2_{H^s}d\tau\leq \|\theta(t_0)\|^2_{H^s}+2C\int^t_{t_0}\|\theta(\tau)\|_{\dot{H}^{2-2\alpha}}\||D|^{\alpha}\theta(\tau)\|^2_{H^s}d\tau$$
Let $t_0>0$ such that $$\|\theta(t)\|_{\dot{H}^{2-2\alpha}}< \frac{1}{2C}, \quad t\geq t_0.$$
Then
$$
\|\theta(t)\|^2_{H^s}+ \int^t_{t_0}\||D|^{\alpha}\theta\|^2_{H^s}d\tau\leq \|\theta(t_0)\|^2_{H^s}$$
Let $t_0>0$ such that $$\|\theta(t)\|_{\dot{H}^{2-2\alpha}}< \frac{1}{2C}, \quad t\geq t_0.$$
Which implies $\mathcal{C} (\R^+, H^s(\mathbb{R}^2))\cap L^2(\R^+, \dot{H}^{s+\alpha}(\mathbb{R}^2)).$
Let $k_0\in \N, \quad k_0\geq2$, such that $$2-2\alpha+k_0\alpha\leq s< 2-2\alpha+(k_0+1)\alpha.$$
\textbf{First case:} If $s=2-2\alpha+k_0\alpha$: we want to prove that
\begin{align*}
   \begin{cases}
   &\underset{t\rightarrow\infty}{\lim}\|\theta(t)\|_{\dot{H}^{2-2\alpha+i\alpha}}=0,\\
   &\int^\infty_0\||D|^\alpha\theta\|^2_{\dot{H}^{2-2\alpha+i\alpha}}<\infty, \ \ \forall 0\leq i\leq k_0.
   \end{cases}
   \end{align*}
   Suppose that $\underset{t\rightarrow\infty}{\lim}\|\theta(t)\|_{\dot{H}^{2-2\alpha+i\alpha}}=0,\ 0\leq \forall i\leq k_0-1 $ and we prove that
   $\underset{t\rightarrow\infty}{\lim}\|\theta(t)\|_{\dot{H}^{2-2\alpha+(i+1)\alpha}}=0$. Let $t_0>0$ such that $$\|\theta(t)\|_{\dot{H}^{2-2\alpha}}< \frac{1}{4C_\alpha}, \quad t\geq t_0.$$
   As $\int^\infty_0\||D|^\alpha\theta\|^2_{\dot{H}^{2-2\alpha+i\alpha}}<\infty$ we have $$\int^\infty_{t_0}\||D|^\alpha\theta\|^2_{\dot{H}^{2-2\alpha+i\alpha}}<\infty,$$
   which implies that there is $t_1>t_0$ such that $$\|\theta(t_1)\|^2_{\dot{H}^{2-2\alpha+(i+1)\alpha}}<\varepsilon.$$
   For $t\geq t_1$ we get
    $$
\|\theta(t)\|^2_{H^{2-2\alpha+(i+1)\alpha}}+2 \int^t_{t_1}\||D|^{\alpha}\theta\|^2_{H^{2-2\alpha+(i+1)\alpha}}\leq \|\theta(t_0)\|^2_{H^s}+2C_\alpha\int^t_{t_1}\|\theta\|_{\dot{H}^{2-2\alpha}}\||D|^{\alpha}\theta\|^2_{H^{2-2\alpha+(i+1)\alpha}}.$$
   Then $$
\|\theta(t)\|^2_{H^{2-2\alpha+(i+1)\alpha}}+ \int^t_{t_1}\||D|^{\alpha}\theta\|^2_{H^{2-2\alpha+(i+1)\alpha}}\leq \|\theta(t_1)\|^2_{H^{2-2\alpha+(i+1)\alpha}}<\varepsilon^2, \quad \forall t\geq t_1.$$
Thus, we get for all $k_0>2$,
\begin{align*}
   \begin{cases}
   &\lim_{t\rightarrow\infty}\|\theta(t)\|_{\dot{H}^{2-2\alpha+k_0\alpha}}=0,\\
   &\int^\infty_0\||D|^\alpha\theta\|^2_{\dot{H}^{2-2\alpha+k_0\alpha}}<\infty.
   \end{cases}
   \end{align*}
   \textbf{Second case:} If $2-2\alpha+k_0\alpha< s< 2-2\alpha+(k_0+1)\alpha$
   Let $\varepsilon>0$ and $r_0\geq0$ such that $$\|\theta(t)\|_{\dot{H}^{2-2\alpha}}< \min(\frac{1}{4C_s}, \varepsilon/2), \quad t\geq r_0.$$
   For $t>a>r_0$, we have $$
\|\theta(t)\|^2_{H^s}+2 \int^t_{a}\||D|^{\alpha}\theta\|^2_{H^s}d\tau\leq \|\theta(t_0)\|^2_{H^s}+2C\int^t_{a}\|\theta(\tau)\|_{\dot{H}^{2-2\alpha}}\||D|^{\alpha}\theta(\tau)\|^2_{H^s}d\tau$$
Then
$$
\|\theta(t)\|^2_{H^s}+ \int^t_{a}\||D|^{\alpha}\theta\|^2_{H^s}d\tau\leq \|\theta(a)\|^2_{H^s}.$$
It suffice to prove that, there is $a\geq t_0$ such that $$\|\theta(a)\|^2_{H^s}\leq \varepsilon.$$
By interpolation, we obtain
$$\|\theta(t)\|_{H^{s}}\leq \|\theta(t)\|^{1-r}_{H^{2-2\alpha+k_0\alpha}}\|\theta(t)\|^{r}_{H^{2-2\alpha+(k_0+1)\alpha}},$$
with $r \in [0, 1[$, which implies
$$\|\theta(t)\|^{\frac{2}{1-r}}_{H^{s}}\leq \|\theta(t)\|^{\frac{2r}{1-r}}_{H^{2-2\alpha+k_0\alpha}}\|\theta(t)\|^{2}_{H^{2-2\alpha+(k_0+1)\alpha}}.$$
Let $\varepsilon >0$
and $E_{\varepsilon }= \{t\geq r_0, \|\theta(t)\|_{\dot{H}^s}>\varepsilon\}$ we get
$$M_0\geq\int_{E_{\varepsilon }}\|\theta(t)\|^{\frac{2}{1-r}}_{\dot{H}^s}dt\geq \varepsilon^{\frac{2}{1-r}} \lambda_1(E_{\varepsilon }).$$
Then $\lambda_1(E_{\varepsilon })\leq t_{\varepsilon}=\varepsilon^{\frac{1-r}{2}}M_0$. For $\eta>0$, there exists $ a\in [t_0,  t_{\varepsilon}+\eta] $ such that $a \notin E_{\varepsilon }$ and it results that
$$\|\theta(a)\|_{\dot{H}^s}\leq\varepsilon,$$
  which prove the desired result.\\
   %%%%%%%%%%%%%%%%%%%%%%%%%%%%%%%%%%%%%%%%
{\bf Acknowledgements.}
It is pleasure to thank J. Benameur for insightful comments and assistance through this work.

\end{document}